\documentclass[a4paper,11pt,reqno]{amsart} 
\pdfoutput=1

\usepackage{fullpage}

\usepackage[T1]{fontenc}
\usepackage[utf8]{inputenc}

\usepackage{lmodern}

\usepackage{microtype} 

\usepackage{url} 

\usepackage[version=3]{mhchem} 
\mhchemoptions{textfontcommand=\sffamily}
\mhchemoptions{mathfontcommand=\mathsf}

\usepackage{amsmath}
\usepackage{amssymb}

\usepackage{graphicx}
\usepackage{caption}
\usepackage{subcaption}
\usepackage{float}
\usepackage{wrapfig}

\usepackage{multirow}
\usepackage{array} 

\usepackage{mathrsfs} 

\usepackage{bm}

\usepackage{bbm} 

\usepackage{enumitem}

\usepackage[dvipsnames]{xcolor}

\usepackage[perpage]{footmisc}

\usepackage{tikz-cd}

\usepackage[all]{xy}

\usepackage[unicode,psdextra,pdfusetitle,
bookmarks, bookmarksdepth=2, colorlinks=true, linkcolor=blue!50!black, citecolor=Black!70, urlcolor=blue!50!black]{hyperref}

\usepackage{bookmark}
\hypersetup{bookmarksnumbered}

\usepackage[indent,skip=.2\baselineskip]{parskip}
\usepackage[nameinlink]{cleveref}
\crefname{ex}{Example}{Examples}
\crefname{thm}{Theorem}{Theorems} 
\crefname{lem}{Lemma}{Lemmas}
\crefname{prop}{Proposition}{Propositions}
\crefname{cor}{Corollary}{Corollaries} 
\crefname{con}{Conjecture}{Conjectures} 
\crefname{def}{Definition}{Definitions}
\crefname{rmk}{Remark}{Remarks}
\crefname{thmalph}{Theorem}{Theorems}

\usepackage{amsthm}

\newtheorem{theorem}{Theorem}[section]
\newtheorem{lemma}[theorem]{Lemma}
\newtheorem{proposition}[theorem]{Proposition}
\newtheorem{corollary}[theorem]{Corollary}

\newtheorem{theoremalphabetic}{Theorem}

\theoremstyle{definition}
\newtheorem{definition}[theorem]{Definition}
\newtheorem{example}[theorem]{Example}

\newtheorem{remark}[theorem]{Remark}

\newtheorem*{notation*}{Notation}

\numberwithin{equation}{section}

\newcommand{\from}{\ensuremath{\colon}}

\DeclareMathOperator{\rk}{rk} 
\DeclareMathOperator{\im}{im} 
\DeclareMathOperator{\spn}{span} 
\DeclareMathOperator{\id}{id} 
\DeclareMathOperator{\diag}{diag} 
\DeclareMathOperator{\Lin}{Lin} 
\DeclareMathOperator{\rad}{rad} 

\renewcommand{\S}{\mathcal{S}}
\newcommand{\B}{\mathcal{B}}

\newcommand{\D}{\mathcal{D}}

\newcommand{\F}{\mathcal{F}}
\newcommand{\X}{\mathcal{X}}

\newcommand{\A}{\mathcal{A}}
\renewcommand{\P}{\mathcal{P}}
\newcommand{\Z}{\mathcal{Z}}

\newcommand{\ZZ}{\ensuremath{\mathbb{Z}}}
\newcommand{\QQ}{\ensuremath{\mathbb{Q}}}
\newcommand{\RR}{\ensuremath{\mathbb{R}}}
\newcommand{\CC}{\ensuremath{\mathbb{C}}}

\newcommand{\I}{\ensuremath{\mathcal{I}}}

\newcommand{\VV}{\mathbb{V}}

\DeclareTextFontCommand{\bfemph}{\bfseries\em}
\newcommand{\term}{\bfemph}

\begin{document}

\title[Generic consistency and nondegeneracy of vertically parametrized systems]{Generic consistency and nondegeneracy of\\ vertically parametrized systems}

\author{Elisenda Feliu, Oskar Henriksson, and  Beatriz Pascual-Escudero}

\begin{abstract}
We determine the generic consistency, dimension and nondegeneracy of the zero locus over $\CC^*$, $\RR^*$ and $\RR_{>0}$ of vertically parametrized systems: parametric polynomial systems consisting of linear combinations of monomials scaled by free parameters. These systems generalize sparse systems with fixed monomial support and freely varying parametric coefficients. As our main result, we establish the equivalence among three key properties: the existence of nondegenerate zeros, the zero set having generically the expected dimension, and the system being generically consistent. Importantly, we prove that checking whether a vertically parametrized system has these properties amounts to an easily computed matrix rank condition.

\vspace{-1.5em}

\end{abstract}

\vspace*{-0.5em}

\maketitle

\section{Introduction}

Polynomial systems that arise in applications often have fixed support, with coefficients that depend on unknown parameters, and a common theme in applied algebraic geometry is to ask how the geometry of the solution set varies with the value of the parameters. 
A particularly well-studied setting is when all coefficients vary independently from each other. We will refer to such systems as \term{freely parametrized systems}. A simple example of this is
\begin{equation}\label{eq:free_example}\left(\begin{array}{l} 
a_{1} x_{1}^2 x_{2}^2 + a_{2} x_{1}^2 x_{3} + a_{3} x_{3} x_{4}^2\\[0.2em]
 a_{4} x_{1}^2 x_{3} + a_{5} x_{3} x_{4}^2 + a_{6} x_{1} x_{3} x_{4} \end{array}\right)\end{equation}
with parameters $a=(a_1,\ldots,a_6)$ and variables $x=(x_1,x_2,x_3,x_4)$.
Many results are known about how the geometry of zero sets of freely parametrized systems for generic choices of parameters depends on the  geometry of the associated Newton polytopes. For instance, Bernstein's theorem describes the generic number of roots in $(\CC^*)^n$ in the square case \cite{bernstein}, and more recent work characterizes generic nonemptiness and generic irreducibility of the zero locus of (possibly non-square) freely parametrized systems \cite{Yu2016prime,khovanskii2016newton}. 

In many applications, the systems have additional structure apart from a fixed support, which gives rise to algebraic dependencies between the coefficients (in particular, some coefficients might be fixed). This happens for instance in enumerative geometry \cite{EH16}, optimization \cite{lindberg2023algebraic}, statistics \cite{HS14}, dynamics \cite{breiding2022TheAD}, reaction network theory~\cite{Obatake-Newton}, and for extremal-generic systems \cite{spaenlehauer-bender}.

In this work, we will study a particular generalization of freely parametrized systems that we refer to as \term{vertically parametrized systems} (a terminology that first appeared in \cite{helminck2022generic}), where we allow linear dependencies among the coefficients of terms with the same monomial. More specifically, a vertically parametrized system is one that can be written as
\[F=C(a\star x^M)\in \CC[a_1,\ldots,a_m,x_1^\pm,\ldots,x_n^\pm]^s\]
where $M\in\ZZ^{n\times m}$ is a matrix encoding monomials, the parameters $a=(a_1,\ldots,a_m)$ scale the monomials via component-wise multiplication $\star$, and $C\in\RR^{s\times m}$ is a matrix of full row rank, whose rows encode linear combinations of the scaled monomials.
A simple example of a vertically parametrized system with the same support as \eqref{eq:free_example} is
\begin{equation}\label{eq:vertical_example}\left(\begin{array}{l}
a_{1} x_{1}^2 x_{2}^2 + 3 a_{2} x_{1}^2 x_{3} + a_{3} x_{3} x_{4}^2\\[0.2em]
 a_{2} x_{1}^2 x_{3} + 2 a_{3} x_{3} x_{4}^2 + a_{4} x_{1} x_{3} x_{4}
\end{array}\right).
\end{equation}

More generally, we will be interested in \emph{linear sections with fixed direction} of the zeros of a vertically parametrized system, in the sense that we also include $\ell\geq 0$ affine forms $Lx-b$ for a fixed matrix $L\in\CC^{\ell\times n}$, and parameters $b=(b_1,\ldots,b_\ell)$, to obtain an \term{augmented vertically parametrized system} of the form
\[F=\left(C(a\star x^M),\, Lx-b\right)\in\CC[a_1,\ldots,a_m,b_1,\ldots,b_\ell,x_1^\pm,\ldots,x_n^\pm]^{s+\ell}.\]
For instance, an example of such a linear section of \eqref{eq:vertical_example} is given by
\begin{equation}\label{eq:augmented_vertical_example}\left(\begin{array}{l} a_{1} x_{1}^2 x_{2}^2 + 3 a_{2} x_{1}^2 x_{3} + a_{3} x_{3} x_{4}^2\\[0.1em]
 a_{2} x_{1}^2 x_{3} + 2 a_{3} x_{3} x_{4}^2 + a_{4} x_{1} x_{3} x_{4}\\[0.1em]
x_2+2x_3-b_1\\[0.1em]
x_2+x_3-b_2
\end{array}\right).
\end{equation}

One of the first main results of this work is a characterization of when an augmented vertically parametrized system is \emph{generically consistent}, in the sense that the variety $\VV_{\CC^*}(F_{a,b})\subseteq(\CC^*)^n$ is nonempty for generic choices of parameters $(a,b)\in\CC^m\times\CC^\ell$. The key condition is the following:
\begin{equation}\label{eq:rank_condition}
        \rk\Big(\begin{bmatrix} C\diag(w)M^\top\diag(h)\\ L\end{bmatrix}\Big)=s+\ell\quad\text{for some $(w,h)\in\ker(C)\times(\CC^*)^n$}.
\end{equation}
This condition ensures that $F_{a,b}$ has a nondegenerate zero (i.e., a zero where the Jacobian of $F_{a,b}$ does not have full rank) for some $(a,b)\in\CC^m\times\CC^\ell$. Using this, we will show in \Cref{ex:examples_from_intro_revisited} that the systems \eqref{eq:free_example} and \eqref{eq:vertical_example} are generically consistent, whereas \eqref{eq:augmented_vertical_example} is generically inconsistent.

\begin{theoremalphabetic}[\Cref{thm:main_theorem_for_vertical_systems}]\label{thmalph:main_result_over_C}
    Let $F=(\,C(a\star x^M),\, Lx-b\,)$ be an augmented vertically parametrized system with $C\in\CC^{s\times m}$ of full row rank, $M\in\ZZ^{n\times m}$ and $L\in\CC^{\ell \times n}$. Then there   are precisely two possibilities: 
    \begin{enumerate}[label=(\roman*)]
        \item \textbf{Generic consistency:} If \eqref{eq:rank_condition} holds, then for generic $(a,b)\in\CC^m\times\CC^\ell$, it holds that $\VV_{\CC^*}(F_{a,b})$ is nonempty of pure dimension $n-(s+\ell)$, and all zeros are nondegenerate.
        \item \textbf{Generic inconsistency:} If \eqref{eq:rank_condition} does not hold, then $\VV_{\CC^*}(F_{a,b})$ is empty for generic $(a,b)\in\CC^m\times\CC^\ell$. Whenever $\VV_{\CC^*}(F_{a,b})$ is nonempty, its dimension is strictly greater than $n-(s+\ell)$ and all zeros are degenerate.
    \end{enumerate}
    Furthermore, the ideal $\langle F_{a,b}\rangle\subseteq\CC[x_1^\pm,\ldots,x_n^\pm]$ is radical for generic $(a,b)\in\CC^m\times\CC^\ell$.
\end{theoremalphabetic}
 
For a freely parametrized system with supports $\mathcal{S}_1,\ldots,\mathcal{S}_s\subseteq\ZZ^n$, we see in \Cref{cor:freely_parametrized} that \eqref{eq:rank_condition} specializes to the  following condition, which has previously been proven by other means in \cite[Lemma~1]{Yu2016prime} and \cite[Theorem~11]{khovanskii2016newton}:
\begin{equation}\label{eq:free_condition}
    \text{There exists a linearly independent tuple $(u_1,\ldots,u_s)\in\textstyle\prod_{i=1}^s\Lin(\mathcal{S}_i$)}\,.
\end{equation}
Here, $\Lin(\S_i)\subseteq\RR^n$ denotes the direction of the affine hull of $\S_i$.
In the square case $s=n$, this, together with Bernstein's theorem, also recovers the usual characterization of when the mixed volume is nonzero (see, e.g. \cite[Theorem~5.1.7]{schneider2013convex}).

In many applications, it is natural to require the parameters or the variables to be (positive) real numbers. If a vertically parametrized system has been defined over the real numbers, in the sense that $C$ and $L$ have real entries, we have the following real version of \Cref{thmalph:main_result_over_C} (a more general result that encompasses both these versions is given in \Cref{thm:main_theorem_for_vertical_systems}).

\begin{theoremalphabetic}[\Cref{thm:main_theorem_for_vertical_systems}]\label{thmalph:main_result_over_R}
    Let $F=(\,C(a\star x^M),\, Lx-b\,)$ be an augmented vertically parametrized system, with $C\in\RR^{s\times m}$ of full row rank, $M\in\ZZ^{n\times m}$ and $L\in\RR^{\ell \times n}$. Suppose that $\ker(C)\cap\RR^m_{>0}\neq\varnothing$. Then there are two cases:
    \begin{enumerate}
        \item[(i)] \textbf{Consistency in a Euclidean open set:} If condition \eqref{eq:rank_condition} holds, then $\VV_{\CC^*}(F_{a,b})\cap\RR^n_{>0}\neq\varnothing$ for $(a,b)$ in a nonempty Euclidean open subset of $\RR^m_{>0}\times\RR^\ell$. For generic such parameter values, all zeros are nondegenerate and it holds that $\dim(\VV_{\CC^*}(F_{a,b})\cap\RR^n_{>0})=n-(s+\ell)$ as a semialgebraic set.
        \item[(ii)] \textbf{Generic inconsistency:} If condition \eqref{eq:rank_condition} does not hold, then the set of $(a,b)\in\RR^m_{>0}\times\RR^\ell$ for which $\VV_{\CC^*}(F_{a,b})\cap\RR^n_{>0}\neq\varnothing$ is nonempty but has empty Euclidean interior, and all zeros are degenerate.
    \end{enumerate}
\end{theoremalphabetic}

By specializing \Cref{thmalph:main_result_over_R} to real freely parametrized systems with prescribed signs for the coefficients, we find in \Cref{cor:fixedsignfree} 
that the zero locus contains positive points (i.e., in $\RR^n_{>0}$) for coefficients in a nonempty Euclidean open set exactly when \eqref{eq:free_condition} holds and each polynomial has at least one positive and one negative coefficient.

There are two key features of augmented vertically parametrized systems that we use in the proof of \Cref{thm:main_theorem_for_vertical_systems}. The first of them is the well-behaved geometry of the incidence variety 
\[\I:=\Big\{(a,b,x)\in\CC^m\times\CC^\ell\times(\CC^*)^n:C(a\star x^M)=Lx-b=0\Big\},\]
which allows us to use classical results from algebraic geometry about the dimension of fibers to relate generic consistency of the system with the generic dimension of the solution set.

\begin{theoremalphabetic}[\Cref{thm:irreducible_smooth_incidence_variety}]\label{thmalph:irreducibility_of_incidence_variety}
For an augmented vertically parametrized system, the incidence variety $\I$ admits a rational parametrization, and is a nonsingular irreducible variety of dimension $m+n-s$. 
\end{theoremalphabetic}

This generalizes the well-known fact that the incidence variety is irreducible for freely parametrized systems (see, e.g., \cite[Proof of Proposition~2.3]{Pedersen-Sturmfels}).   In the square case (when $s+\ell=n$), another important consequence of \Cref{thmalph:irreducibility_of_incidence_variety} is that the  Galois action of the monodromy group is transitive, which ensures that augmented vertically parametrized systems can be solved numerically with monodromy methods \cite{Duff2018monodromy}. 

As a second key ingredient,   we show in \Cref{lem:critical_points} that degenerate zeros correspond to the critical points of the parametrization of $\I$. This allows us to invoke Sard's lemma and relate the dimension of the complex and real varieties, and thereby derive \Cref{thmalph:main_result_over_R} from \Cref{thmalph:main_result_over_C}. 

Although  \Cref{thmalph:main_result_over_C} and \Cref{thmalph:irreducibility_of_incidence_variety} are stated for zeros in the complex torus $(\CC^*)^n$, we show in \Cref{thm:extension_to_C} that if $M\in\ZZ^{n\times m}_{\geq 0}$ and each polynomial contains an independent constant term, then both results also hold for the full zero locus in $\CC^n$. In fact, we prove that an augmented vertically parametrized system with independent constant terms generically does not have any irreducible component contained in a coordinate hyperplane. This generalizes  \cite[Lemma~2.1]{LW96} from the freely parametrized and square case to the vertically parametrized setting.
 
\smallskip
The conclusions of this work can be contrasted with the following non-vertically parametrized systems, where 
neither case (i) nor (ii) of \Cref{thmalph:main_result_over_C} apply:
\begin{enumerate}[label=(\roman*)]
    \item The following system is \emph{horizontally parametrized} in the language of \cite{helminck2022generic} (in the sense that we allow linear dependencies among coefficients appearing within the same equation):
    \begin{equation}\label{eq:non_vertical_systems1}
    \left(\begin{array}{l} a_{1} x_{1}^{2}-a_{1} x_{1}-a_{3} x_{1}+a_{3} \\ a_{2} x_{1} x_{2}-a_{2} x_{2}-a_{4} x_{1}+a_{4} \end{array}\right)=\left(\begin{array}{l} (x_1-1)(a_1x_1-a_3) \\  (x_1-1)(a_2x_2-a_4) \end{array}\right).
    \end{equation}
  The system is generically consistent but the zero sets have dimension one as they all contain the line $\{x_1=1\}$ (and the nondegenerate point $(a_3/a_1,a_4/a_2)$ whenever $a_1a_2\neq 0$).
    \item In the system 
    \begin{equation}\label{eq:non_vertical_systems2}
    \left(\begin{array}{l}a_{1} x_{1}^{2}+a_{3} x_{2}+a_{4} x_{2} x_{3} \\ 2 a_{1} x_{1} x_{2}+a_{2} x_{1}^{2}+a_{3} x_{2} \\ a_{1} x_{2}^{2}+a_{2} x_{1}^{2}-a_{4} x_{2} x_{3}  \end{array}\right),
    \end{equation}
    the parameter $a_1$ accompanies different monomials in different equations, namely $x_1^2,x_1x_2,x_2^2$. The system is generically consistent with zero-dimensional zero set, but all zeros are degenerate.
    \item Similarly to (i), the second polynomial factors in the  system  \begin{equation}\label{eq:example_nonlinear_coefficients}
    \left(\begin{array}{l} a_1 x_1 - a_2 x_2 \\ a_1^2 x_1^2 - a_2^2 x_2^2\end{array}\right),
    \end{equation} 
    from which it follows that the set of zeros is generically one-dimensional, and that all zeros are degenerate. This system is not linear in the parameters.

\end{enumerate}
The incidence varieties of   systems \eqref{eq:non_vertical_systems1} and \eqref{eq:non_vertical_systems2} are reducible, while it is irreducible for \eqref{eq:example_nonlinear_coefficients}, but of dimension $3>m+n-s$. 

\smallskip
We end the introduction by briefly describing two applications where vertically parametrized systems appear naturally: optimization and dynamical systems.

\subsubsection*{Critical points in optimization:} A freely parametrized single Laurent polynomial can be written as $f=a^\top x^M$ for a full row rank matrix $M\in\ZZ^{n\times m}$ and parameters $a=(a_1,\ldots,a_m)$. Its critical points in $(\CC^*)^n$ are the zeros of the square vertically parametrized system
\begin{equation}\label{eq:critical_system}
    x\star \nabla f =M(a\star x^M)\in\CC[a_1,\ldots,a_m,x_1^\pm,\ldots,x_n^\pm]. 
\end{equation}
It follows from \Cref{thmalph:main_result_over_C} that $f_a$ has critical points in $(\CC^*)^n$ for generic $a\in\CC^m$ if and only if
\begin{equation}\label{eq:rank_condition_critical_points}
    \rk(M\diag(w)M^\top)=n\:\:\:\text{for some $w\in\ker(M)$}.
\end{equation}
If this is the case, the number of critical points is generically finite (c.f. \Cref{ex:critical_points_revisited}).

\subsubsection*{Steady states of reaction networks:} A common model for the dynamics of reaction networks is (generalized) \emph{mass-action kinetics}, where the evolution of some quantities $x$ vary according to an autonomous system of ordinary differential equations of the form 
$$\dot{x}(t)=C(a\star x(t)^M),$$ 
where $C\in\ZZ^{n\times m}$ is called the stoichiometric matrix, $a\in\RR_{>0}^m$ the vector of reaction rate constants, and $M\in\ZZ^{n\times m}$   the kinetic matrix. The steady states are then given by the zeros of an augmented vertically parametrized system
\[\left(C(a\star x^M),\, Lx-b\right)\]
where the rows of $L$ encode conserved quantities (first integrals), and $b$ encodes total amounts. For an introduction to the theory of reaction networks and mass-action kinetics, we refer   to \cite{dickenstein2016inivitation,feinberg2019foundations}, and to \cite{MR12} for the notion of generalized mass-action kinetics.

\begin{example}
The following network from \cite{HMM11} models the interactions between the HIV virus and the T-cells and macrophages of an infected individual:
\[\arraycolsep=1.5em
\begin{array}{llll}
\ce{T + V -> 2T + V} & \ce{M + V -> M_i} & \ce{0 <=> T} & \ce{T_i -> 0 }\\
\ce{T + V -> T_i} & \ce{T_i -> V + T_i} & \ce{0 -> M} & \ce{M -> 0 }\\
\ce{M + V -> 2M + V} & \ce{M_i -> V + M_i} & \ce{V -> 0} &  \ce{M_i -> 0}\,. \\
\end{array}\]
Modeled with mass-action kinetics, we get a stoichiometric and kinetic matrix (the system lacks conserved quantities) that satisfy the conditions in \Cref{thmalph:main_result_over_R}(i). From this, we conclude that there is a nonempty Euclidean open subset of parameter space where the system has a positive steady state. The existence of positive steady states is one of the key features of the model discussed in \cite{HMM11}. If the processes \ce{0 \to T} and \ce{0 \to M} (which correspond to regeneration of T-cells and macrophages) are removed, \Cref{thmalph:main_result_over_R}(ii) applies, and hence the set of parameter values for which there are positive steady states is nonempty but has empty Euclidean interior.
\end{example}

In a follow-up paper {\cite{FeliuHenrikssonPascualCRN}}, we use the results from this paper  to {study the generic geometry of steady state varieties, and }strengthen several previous statements about reaction networks, in particular concerning absolute concentration robustness from {\cite{pascualescudero2020local,ACR_AlgGeom}} and nondegenerate multistationarity from {\cite{Conradi2017identifying}}.

\medskip

The paper is organized as follows. The study of augmented vertically parametrized systems is the content of \Cref{sec:vertical_systems}. Before that, we devote \Cref{sec:parametric_systems} to a more general discussion of parametric polynomial systems with irreducible incidence varieties and where we allow for restricted domains of parameters and variables. We focus on the connection between nondegeneracy of  zeros, generic consistency, the generic dimension of the varieties, and the radicality of the ideals, with the main results being gathered in \Cref{thm:summary_of_implications}.

\subsection*{Notation and conventions}
We let $\star\from\CC^n\times\CC^n\to\CC^n$ denote the Hadamard product, given by $({x}\star{y})_i=x_iy_i$. For a field $K$, we let $K^*=K \setminus \{{0}\}$ be the group of units. For a matrix $A=(a_{ij}) \in\ZZ^{n\times m}$ and a vector ${x}\in(K^*)^n$, we let ${x}^A\in K^m$ be defined by $({x}^A)_j=x_1^{a_{1j}}\cdots x_n^{a_{nj}}$ for $j = 1,\ldots,m$. The zero matrix of size $n\times m$ is denoted $0_{n\times m}$, and ${\rm Id}_n$ is the identity matrix of size $n$.
When saying that a property holds \emph{generically} in a family indexed by some parameters in a set $\P\subseteq \CC^k$, we mean that it holds in a {nonempty open subset of $\P$, with respect to the subspace topology induced by the Zariski topology on $\CC^k$}. For a set $S\subseteq \CC^n$, we let 
$\overline{S}$ denote the Zariski closure of $S$ in $\CC^n$.

\subsection*{Acknowledgements}
EF and OH have been funded by the Novo Nordisk Foundation project with grant reference number NNF20OC0065582. 
BP has been funded by the European Union's Horizon 2020 research and innovation programme under the Marie Sklodowska-Curie IF grant agreement No 794627 and the Spanish Ministry of Economy project with reference number PGC2018-095392-B-I00. This work has  also been funded by the European Union under the Grant Agreement number 101044561,
POSALG. Views and opinions expressed are those of the authors only and do not necessarily
reflect those of the European Union or European Research Council (ERC). Neither the European
Union nor ERC can be held responsible for them. 

The authors thank  Ignacio González Mantecón and Anne Shiu for helpful feedback on earlier versions of the manuscript, and Bernd Sturmfels for pointing us to \cite{Yu2016prime}.

\section{Preliminaries of parametric polynomial systems}\label{sec:parametric_systems}

In this section, we study the generic properties of the zero set of parametric systems, under the assumption that the incidence variety is irreducible of known dimension. This sets the background theory of vertically parametrized systems that we will develop in \Cref{sec:vertical_systems}.

We focus mainly on zeros in $(\CC^*)^n$, $(\RR^*)^n$ and $\RR^n_{>0}$, and we therefore allow Laurent polynomials with negative exponents. However, we will show later that the theory in this section also extends to $\CC^n$ if one restricts to polynomials with nonnegative exponents (see \Cref{rk:extension_C} and \Cref{subsec:extension_to_C}).

\subsection{Framework}\label{subsec:framework}
We consider a parametric  (Laurent) polynomial system of the form
\[
F = (f_1,\dots,f_r) \in\CC[p_1,\dots,p_k, x_1^\pm,\dots,x_n^\pm]^r,
\]
for some integers $k,n,r>0$ with $r\leq n$, where we view $p=(p_1,\ldots,p_{k})$ as parameters in $\CC^k$ and $x=(x_1,\ldots,x_n)$ as variables in $(\CC^*)^n$. Throughout this section, when writing $\CC[p,{x}^{\pm}]$, we implicitly assume that $p$ has $k$ entries and $x$ has $n$ entries. We consider the \term{incidence variety}
\[\I :=\{(p,{x})\in \CC^k\times (\CC^*)^n : F(p,{x})={0}\}\]
and the projection map to parameter space
\[
\pi\from \I \rightarrow \CC^k,\quad  (p,{x})\mapsto p,  
\]
where the system $F$ is implicit in the notation. For each choice of parameters $p\in\CC^{k}$, we get a specialized system
\[
F_{p}:=F(p,\cdot)\in \CC[{x}^{\pm}]^r.
\]
We identify the very affine variety 
\[\VV_{\mathbb{C}^*}(F_{p})=\{x\in(\CC^*)^n:F_{p}(x)=0\}\subseteq(\CC^*)^n\] 
with the fiber $\pi^{-1}(p)$ of the projection map. We also form the set of parameter values for which the system is consistent:
\begin{equation}\label{eq:Z}
\Z:=\{p\in\CC^k:\VV_{\CC^*}(F_p)\neq\varnothing\}=\pi(\I)\subseteq\CC^k.
\end{equation}

{We say that $F$ is \term{generically consistent} if $\Z$ is Zariski dense in $\CC^k$. Since $\Z$ is constructible (see, e.g., \cite[Theorem~3.2.3]{cox2015ideals}), this is equivalent to $\Z$ containing a Zariski open subset of~$\CC^k$. 
Hence, if $F$ is generically consistent, a property holds generically in $\Z$ if and only if it holds generically in $\CC^k$.}

We will restrict the parameter space to a subset $\P\subseteq \CC^k$ (with $\P=\RR^k_{>0}$ as the main example), and we therefore extend the notation in \eqref{eq:Z} to
\begin{equation}\label{eq:Z_P}
\Z_{\P}: = \Z\cap \P =\{ p\in  \P : \VV_{\CC^*}(F_{p})  \neq \varnothing\}=\pi\big(\I\cap (\P\times (\CC^*)^n)\big) \subseteq \CC^k.
\end{equation}
{Under mild assumptions, the choice of $\P$ does not affect generic consistency of the system.}

{
\begin{lemma}
\label{lem:denseness_of_ZP_and_Z}

Let $F\in\CC[p,x^\pm]^r$ and $\P\subseteq\CC^k$. Then the following holds:
\begin{enumerate}
    \item[(i)] If $\overline{\Z_\P}=\CC^k$, then $\overline{\Z}=\CC^k$.
    \item[(ii)] If $\overline{\P}=\CC^k$, then the converse of (i) holds.
\end{enumerate}
\end{lemma}
\begin{proof}
Part (i) is immediate, since $\Z_\P\subseteq\Z$. For part (ii), note that if $\overline{\Z}=\CC^k$, then $\Z$ contains a nonempty Zariski open subset $U$ of $\CC^k$. The intersection $U\cap\P$ is Zariski dense in $\CC^k$ since it is the intersection of a nonempty Zariski open and a Zariski dense set, and $\overline{\Z_\P}=\CC^k$ follows.
\end{proof}
}

With this notation in place, we proceed now to study the generic dimension of $\VV_{\CC^*}(F_p)$ for $p\in\P$ in relation to whether the system is generically consistent.
We then move on to study zeros in subsets $\X\subseteq(\CC^*)^n$. We also relate the study of nonemptiness and dimension of $\VV_{\CC^*}(F_p)$ to the concept of \emph{nondegeneracy}. The connections among  these properties are summarized in 
\Cref{thm:summary_of_implications} towards the end of this section.

\subsection{Generic consistency and dimension}
An immediate first observation is that for each $p\in\P$,  the principal ideal theorem \cite[Theorem~10.2]{eisenbud1995commutative} gives that all irreducible components of $\VV_{\CC^*}(F_{p})$  have dimension at least $n-r$. In particular,
\begin{equation}\label{eq:lower_dimension_bound}
\dim (\VV_{\CC^*}(F_p)) \geq n-r\qquad \text{for all $p\in\Z$}.
\end{equation}
If \eqref{eq:lower_dimension_bound} holds with equality for a given $p\in\P$, then all irreducible components have dimension $n-r$, and we say that $\VV_{\CC^*}(F_p)$ has \term{pure dimension} $n-r$. In what follows, the bound \eqref{eq:lower_dimension_bound} will be referred to as the \term{expected dimension} of $\VV_{\CC^*}(F_{p})$. Likewise, the expected dimension of $\I$ is $k+n-r$. 

\begin{remark}\label{rem:dim}
As neither of the expected dimensions $\dim(\I)=k+n-r$ and $\dim(\VV_{\CC^*}(F_p))=n-r$ can be attained if the coefficient matrix of $F$ (seen as a polynomial system in $\CC[p,{x}^{\pm}]^r$) fails to have full rank, a natural preprocessing step when studying the dimension of the zero set  is to remove linear dependencies in the entries of $F$.  
\end{remark}

{
\begin{theorem}
\label{thm:generic_dimension}
Let $\P \subseteq \CC^k$ and $F\in \CC[p,{x}^{\pm}]^r$ be such that the incidence variety $\I$ is irreducible. Then:
\begin{enumerate}
\setlength{\itemsep}{0.2em}
\item[(i)] For all $p\in\Z$, it holds that \[\text{$\dim(Y)\geq\dim(\I)-\dim(\overline{\Z})$\:\: for all irreducible components $Y\subseteq \VV_{\CC^*}(F_p)$}\]
with equality for generic $p\in\overline{\Z}$.
\item[(ii)]
If $\overline{\Z_{\P}}=\CC^k$, 
then $\VV_{\CC^*} (F_{p})$ has pure dimension $\dim(\I)-k$ for generic $p\in\Z_{\P}$. 
\item[(iii)] If  $\overline{\Z_{\P}}\subsetneq\CC^k$ and $\overline{\P}=\CC^k$, then all irreducible components of $\VV_{\CC^*} (F_{p})$ have dimension greater than $\dim(\I)-k$ for all $p\in\Z$.
\end{enumerate}
\end{theorem}

\begin{proof}
Part (i) follows by applying the classical theorem of dimension of fibers (see, e.g., \cite[Theorem~3.13, Corollary~3.15]{Mumford}) to the canonical projection map 
\[\pi\colon \I \rightarrow \overline{\Z},\quad (p,x)\mapsto p.\]
Part (ii) follows by noting that (i) and \Cref{lem:denseness_of_ZP_and_Z}(i) together give that there is a nonempty Zariski open subset $U\subseteq\CC^k$ such that $\VV_{\CC^*}(F_p)$ has pure dimension $\dim(\I)-k$ for all $p\in U$. In particular, this holds for $p\in\Z_\P\cap U$, which is nonempty since $\overline{\Z_\P}=\CC^k$.
Finally, part (iii) follows directly from (i) and \Cref{lem:denseness_of_ZP_and_Z}(ii).
\end{proof}
}

Combining part (ii) and (iii) of \Cref{thm:generic_dimension} {when $\overline{\mathcal{P}}=\CC^k$}, we obtain that  $\dim(\VV_{\CC^*}(F_p))=\dim(\I)-k$ for some $p\in\Z_\P$ if and only if $\Z_{\P}$ is Zariski dense in $\CC^{k}$.

\begin{example}
For the system in \eqref{eq:example_nonlinear_coefficients} with $n=k=r=2$, we have $\Z=(\CC^*)^2$ and $\I$ is irreducible of dimension $3$. 
By \Cref{thm:generic_dimension}(ii), it holds that $\dim(\VV_{\CC^*}(F_p))=3-2=1$ for generic $p\in(\CC^*)^2$, which we already noticed in \eqref{eq:example_nonlinear_coefficients}. In this case, the generic dimension differs from the expected dimension $n-r=0$ from \eqref{eq:lower_dimension_bound}.
\end{example}
 
{
We end by connecting consistency and dimension with the notion of flatness. Recall that a morphism of varieties $f\from X\to Y$  is said to be flat at a point $x\in X$ if the induced homomorphism of local rings $\mathscr{O}_{Y,f(x)}\to\mathscr{O}_{X,x}$ is a flat ring homomorphism. 

\begin{proposition}
\label{prop:flatness}
Let $F\in\CC[p,x^\pm]^r$ be a parametric system such that $\I$ is irreducible. Then the following holds:
\begin{enumerate} 
    \item[(i)] If the projection $\pi\from\I\to\CC^k$ is flat at generic points $(p,x)\in\I$, then $\overline{\Z}=\CC^k$. 
    \item[(ii)] If $\I$ is nonsingular, then the converse of (i) holds.
\end{enumerate}
\end{proposition}

\begin{proof}
It follows from \cite[00ON]{stacks} that if $\pi\from\I\to\CC^k$ is flat at a point $(p,x)\in\I$, then the local dimension of $\pi^{-1}(p)$ at $(p,x)$ is $\dim(\I)-k$, and  \cite[00R4]{stacks} gives that the converse is true if $\I$ is nonsingular (and hence in particular Cohen--Macaulay). The desired result now follows from \Cref{thm:generic_dimension}.
\end{proof}
}

\subsection{Restricting the ambient space}
We consider now the zero set obtained by restricting the ambient space to a subset $\X\subseteq(\CC^*)^n$. The main example we have in mind is the positive orthant $\X=\RR^n_{>0}$.  
For $K\in\{\RR,\CC\}$, $F\in K[x^\pm]^r$, and $\X\subseteq(K^*)^n$, 
we let the variety $\VV_{K^*}^{\X}(F)$ be defined as the union of the irreducible components of  $\VV_{K^*}(F)$ that intersect $\X$. 
Clearly, the expected dimension $n-r$ from \eqref{eq:lower_dimension_bound} is a lower bound on the dimension of $\VV_{\CC^*}^{\X}(F)$ as long as 
$\VV_{\CC^*}(F)\cap\X\neq\varnothing$.

\begin{remark}\label{rem:density_in_restriction_variety}   
As  $\VV^{\X}_{\CC^*}(F)\subseteq (\CC^*)^n$ is a Zariski closed set, we have
\begin{equation}\label{eq:zarinclusion}
\overline{\VV_{\CC^*} (F)\cap \X} \subseteq \VV^{\X}_{\CC^*}(F),
\end{equation}
but equality might not hold (consider $F=(x_1 - 1)^2 + (x_2 - 1)^2$ and $\X=\RR^2_{>0}$). 
In general,  equality in \eqref{eq:zarinclusion} holds 
if $\X$ is a Euclidean open subset of $(\CC^*)^n$, or  if $\X$ is a Euclidean open subset of $(\RR^*)^n$ and additionally each irreducible component of $\VV_{\CC^*}^{\X}(F)$ contains a nonsingular real point
(see \cite[Theorem~6.5]{pascualescudero2020local} and \cite[Proposition~3.3.16]{BCR}). 
Furthermore, if $F\in\RR[x^\pm]^r$, then $n-r$ is a lower bound of the dimension of $\VV_{\RR^*}(F)$ if it is nonempty and each of its irreducible components contains a nonsingular point. Similarly,  $n-r$ is a lower bound on the dimension of $\VV_{\RR^*}(F)\cap\RR^n_{>0}$ as a semialgebraic set if it is nonempty and each irreducible component of $\VV_{\RR^*}(F)$ that intersects $\RR^n_{>0}$ contains a nonsingular point.
\end{remark} 

For $\P\subseteq\CC^{k}$ and $\X\subseteq(\CC^*)^n$, we generalize definition~\eqref{eq:Z_P} to restrict to zeros in $\X$:
\[\Z_{\P}(\X):= \{ p\in  \P : \VV_{\CC^*}(F_p) \cap \X \neq \varnothing\}=\pi\big(\I\cap (\P\times \X)\big) \subseteq \CC^k,\]
where recall that $\pi$ is the   projection to parameter space.
With this notation, $\Z_{\P}((\CC^*)^n)=\Z_{\P}$, and we have the following inclusions:
\begin{equation}\label{eq:Z_inclusions}
\Z_{\P}(\X) \subseteq  \Z_{\P}\subseteq \Z, \qquad \Z_{\P}(\X) \subseteq  \Z_{\CC^k}(\X)\subseteq \Z.
\end{equation}
If $\Z_{\P}(\X)$ is Zariski dense in $\CC^k$, then so is $\Z_{\P}$. If $\I$ is irreducible and attains its expected dimension $k+n-r$, then \Cref{thm:generic_dimension}(ii) readily gives that $\VV_{\CC^*}^{\X} (F_p)$ has pure dimension $n-r$ for generic $p\in \Z_{\P}(\X)$.  The following example illustrates that \Cref{thm:generic_dimension}(iii) might not extend   to subsets $\X$ unless extra requirements are  imposed.

\begin{example}\label{ex:sum_of_squares}
Let $\P=\RR^3$, $\X=\RR^2_{>0}$, and 
$f = (p_1 x_1-p_2x_2)^2 + p_3^2x_1.$
The incidence variety $\I$ is irreducible and attains its expected dimension $4$, and $\VV_{\CC^*}(f_p)$ has dimension $1$ for generic $p$. The set 
\[\Z_{\P}(\X) = \{ p\in \RR^3 \colon p_3=0, p_1p_2> 0\}\] 
is not Zariski dense  in $\CC^3$, but still, the dimension of 
$\VV_{\CC^*}^{\X}(f_p)$ takes the expected value of 1 for all $p \in \Z_{\P}(\X)$. 
\end{example}

In what follows, the sets $\P$ and $\X$ will  be required to satisfy that  $\I\cap (\P\times \X)$ is Zariski dense in $\I$, which ensures the following property.

\begin{lemma}\label{lem:denseness_of_ZPX}
Let $F\in \CC[p,{x}^{\pm}]^r$, $\P\subseteq\CC^{k}$, and $\X\subseteq (\CC^*)^n$ be
such that $\I\cap (\P\times \X)$ is Zariski dense in $\I$. Then 
\[ \overline{\Z_{\P}(\X)}  = \overline{\Z_{\P}}  = \overline{\Z_{\CC^k}(\X)} = \overline{\Z} .\]
In particular $\Z$ is Zariski dense in $\CC^k$ if and only if $\Z_{\P}(\X)$ is.
\end{lemma}
\begin{proof}
As $\I\cap (\P\times \X)$ is Zariski dense in $\I$, we have 
\[\overline{\Z_{\P}(\X)} = \overline{\pi(\I\cap (\P\times \X))}=\overline{\pi(\overline{\I\cap (\P\times \X)})}=\overline{\pi(\I)}=
\overline{\Z},
\]
which together with the inclusions in \eqref{eq:Z_inclusions} give the desired statement. 
\end{proof}

\begin{remark}\label{rk:PxX_dense_in_I}
Assume that $\I$ is nonsingular and irreducible. If $\P$ and $\X$ are Euclidean open subsets of $\CC^k$ and $(\CC^*)^n$, respectively, then $\I\cap (\P\times \X)$  is Zariski dense in $\I$ if this intersection is nonempty. The same conclusion  holds if $F$ has real coefficients and $\P$ and $\X$ are Euclidean open subsets of $\RR^k$ and $\RR^n$, respectively. This follows from the nonsingularity of $\I$, together with the fact that for a complex irreducible variety defined by polynomials with real coefficients that has at least one real nonsingular point, the real points form a Zariski dense subset of the complex variety \cite[Proposition~3.3.16]{BCR}. Observe that in \Cref{ex:sum_of_squares}, the intersection $\I\cap (\P\times \X)$, which is not dense in $\I$, consists only of singular points.
\end{remark}

\subsection{Nondegeneracy, Euclidean interior, and (real) dimension}\label{subsec:nondegeneracy_and_real_dimension}
This subsection recalls the concept of  \emph{nondegeneracy}, as a means of studying the dimension of a variety, both theoretically and computationally. Additionally, we relate Zariski denseness of $\Z_{\P}(\X)$ in $\CC^{k}$ to $\Z_{\P}(\X)$ having nonempty Euclidean interior in $\P$. 

\begin{definition}
Given a polynomial system $F=(f_1,\ldots,f_r)\in \CC[{x^{\pm}}]^r$,  a zero ${x}^*\in\VV_{\CC^*}(F)$ is called \term{nondegenerate} if the Jacobian matrix $J_{F}({x}^*):=\big(\tfrac{\partial f_i}{\partial x_j}(x^*)\big)_{ij}$ has rank $r$. Otherwise, $x^*$ is called \term{degenerate}.
\end{definition}

Unlike the weaker notion of \emph{nonsingularity}, nondegeneracy depends on the particular tuple $F \in\CC[x^{\pm}]^r$, and not just the variety $\VV_{\CC^*}(F)$. 
The relationship between nondegenerate zeros and  nonsingular points, and the connection to  dimension, is summarized in the following proposition, which gathers several well-known results.

\begin{proposition}\label{prop:propoerties_of_nondegenerate_zeros}
\leavevmode
\begin{enumerate}[label=(\roman*)]
\item[(i)]
Let $F\in \CC[x^{\pm}]^r$ and $x^*\in (\CC^*)^n$ be a nondegenerate zero of $F$. Then $x^*$ is a non\-sin\-gular point of $\VV_{\CC^*}(F)$ and belongs to a unique irreducible component, which has dimension $n-r$.
\item[(ii)]
Let $F\in \RR[x^{\pm}]^r$ and $x^*\in (\RR^*)^n$ be a nondegenerate zero of $F$. Then there is a unique irreducible component of $\VV_{\RR^*}(F)$ containing ${x}^*$. This irreducible component is Zariski dense in the irreducible component of $\VV_{\CC^*}(F)$ containing $x^*$ and has dimension $n-r$.
\item[(iii)] Let $F\in\CC[p,x^{\pm}]^r$. Then the set $\mathcal{U}$ of points $(p,x)\in \I$  for which $x$ is a nondegenerate zero of $F_p$ is a Zariski open subset of $\I$.
\end{enumerate}
\end{proposition}
\begin{proof}
Statement (i) is \cite[Theorem~9.6.9]{cox2015ideals}. For (ii), the uniqueness follows from \cite[Theorem~9.6.8(iv)]{cox2015ideals}, whereas the density and dimensionality claims follow from \Cref{rem:density_in_restriction_variety} (see also \cite[Theorem~3.3.10]{BCR}). For statement (iii), consider the Jacobian $J_{F}=(\partial f_i/\partial x_j)\in\CC[p,x^{\pm}]^{r\times n}$, and note that the complement $\I\setminus\mathcal{U}\subseteq\I$ is cut out by the $r$-minors of $J_F$, and therefore forms a Zariski closed subset of $\I$.
\end{proof}
 
\Cref{thm:generic_dimension} allows us to determine the generic dimension of $\VV_{\CC^*}(F_p)$ from the Zariski denseness 
of $\Z_{\P}$, while the implicit function theorem (as we will use below in \Cref{prop:degdim})
gives that nondegeneracy implies that $\Z_{\P}$ has nonempty Euclidean interior. The key for making the connection between consistency and the Euclidean interior is the following topological property.

\begin{definition}
A nonempty subset $A\subseteq\CC^k$ is said to be \term{locally Zariski dense} (in $\CC^{k}$) if
$U\cap A$ is Zariski dense in $\CC^k$ for any 
  Euclidean open subset $U\subseteq\CC^k$ with $U\cap A\neq\varnothing$.
  \end{definition}

Simple examples of locally Zariski dense sets in $\CC^k$ include all nonempty Euclidean open subsets of $\CC^k$, $\RR^k$, and $\RR^k_{>0}$. In general, any nonempty subset $S\subseteq \RR^k$ for which, with the Euclidean topology, its closure agrees with the closure of its interior, is locally Zariski dense.
Any locally Zariski dense set is in particular Zariski dense, but the converse is not true (for instance, $\mathbb{Z}^k$ is Zariski dense in $\CC^k$, but not locally Zariski dense). 

\begin{lemma}\label{lem:nonemptyinterior}
Let $\P\subseteq\CC^{k}$ be locally Zariski dense. If a subset $S\subseteq \P$ has nonempty Euclidean interior in $\P$, then $S$ is Zariski dense in $\CC^k$.  
\end{lemma}

\begin{proof}
By hypothesis, there exists an open Euclidean ball $B\subseteq \CC^k$ such that $\varnothing \neq B\cap \P \subseteq S$.  The Zariski closures satisfy $\CC^k=\overline{B\cap \P}\subseteq \overline{S}$ as $\P$ is locally Zariski dense. Hence $\CC^k=\overline{S}$.
\end{proof}

\begin{proposition}\label{prop:degdim}
Let $F\in \CC[p,{x}^{\pm}]^r$ and $\P \subseteq \CC^k$ be  locally Zariski dense. Assume that  $F_{p^*}$ has a nondegenerate zero $x^*$ in $(\CC^*)^n$ for some  $p^* \in \P$. Then the following statements hold:
\begin{enumerate}
    \item[(i)] $\Z_{\P}$ has nonempty Euclidean interior in $\P$ and is Zariski dense in $\CC^k$.
 
    \item[(ii)]  If in addition the incidence variety $\I$
is irreducible, 
then $\dim(\I)=k+n-r$, and $\VV_{\CC^*}(F_p)$ has pure dimension $n-r$ for generic $p\in\Z_{\P}$.
\end{enumerate}
\end{proposition}
\begin{proof}
For (i), by assumption we have $\rk(J_{F_{p^*}}(x^*))=r$. Let $A\in\CC^{(n-r)\times n}$ be a matrix whose rows extend the rows of $J_{F_{p^*}}(x^*)$ to a basis of $\CC^n$. Then $x^*$ is a nondegenerate zero of the square system $\widetilde{F}_{p^*}$, where
\[\widetilde{F}:=\begin{pmatrix}F \\ Ax-Ax^*\end{pmatrix}\in \CC[p,x^{\pm}]^n.\]
The complex implicit function theorem \cite[Proposition~1.1.11]{huybrechts2005complex} now gives that there exists an open Euclidean neighborhood $B$ of $p^*$ contained in  $\widetilde{\Z}$ (the set $\Z$ for $\widetilde{F}$) and hence in 
$\Z$. Intersecting $B$ with $\P$, we obtain the first part of (i). By \Cref{lem:nonemptyinterior}, $\Z_{\P}$ is  Zariski dense in $\CC^k$.  

For (ii), as the pair $(p^*,x^*)$ is a nondegenerate zero of $F$ and $\I$ is irreducible, \Cref{prop:propoerties_of_nondegenerate_zeros}(i) gives that $\I$ has the expected dimension. This fact,  (i)
and \Cref{thm:generic_dimension}(ii) give now the second part of (ii).
\end{proof}

In certain scenarios within the setting of \Cref{thm:generic_dimension}, all zeros of $F_p$ will be nondegenerate for generic $p\in\P$. A simple such condition is that the system is square.

\begin{theorem}\label{thm:nondegF}
Let $F\in \CC[p,{x}^{\pm}]^r$  with $n=r$,
and let $\P\subseteq\CC^{k}$ be locally Zariski dense. Assume that the incidence variety $\I$
is irreducible 
and that $F_{p^*}$ has a nondegenerate zero  for some $p^*\in \P$. 
Then all  zeros of $F_p$ are nondegenerate for all $p$ in a nonempty Zariski open subset of $\Z_{\P}$.
\end{theorem}
\begin{proof}
Consider the proper Zariski closed subset set $D\subsetneq\I$ consisting of  the points $(p,x)$ for which $x$ is a degenerate zero of $F_p$ (c.f. \Cref{prop:propoerties_of_nondegenerate_zeros}(iii)). As $\dim(\I)=k$ by \Cref{prop:degdim}(ii), we have  $\dim(\overline{\pi(D)})  \leq \dim(D)<k$.
For all $p$ in the nonempty Zariski open set
$U:=\CC^k\setminus \overline{\pi(D)}$,  all zeros of $F_p$ are nondegenerate. The result now follows from $U\cap \Z_{\P}\neq \varnothing$, as $\Z_{\P}$ is Zariski dense by \Cref{prop:degdim}(i). 
\end{proof}

We conclude this subsection by noting that nondegeneracy 
allows us to also assert radicality of the ideals generated by the systems. The following proposition follows from standard commutative algebra arguments. We give a proof in \Cref{app:radical} for completeness.

\begin{proposition}
\label{prop:generic_radicality_in_Laurent_polynomial_ring}
Let $F\in \CC[p,x^{\pm}]^r$, and suppose that all zeros of $F_p$ in $(\CC^*)^n$ are nondegenerate for generic $p\in \CC^{k}$. Then the following holds:
\begin{enumerate}
\item[(i)] The ideals
$\langle F_p \rangle\subseteq \CC[x^\pm]$ and $\langle F_p \rangle\cap \CC[x]\subseteq\CC[x]$ are radical for generic $p\in\CC^k$.
\item[(ii)] The ideals $\langle F\rangle\subseteq\CC(p)[x^\pm]$ and $\langle F\rangle\cap \CC(p)[x]\subseteq\CC(p)[x]$ are radical.
\end{enumerate}
\end{proposition}

\subsection{Main implications on dimension, Zariski denseness and nondegeneracy}
In the previous subsections, we have studied the generic dimension of $\VV_{\CC^*}(F_p)$ in relation to the existence of nondegenerate zeros and topological properties of $\Z_{\P}(\X)$. 
For $F\in \CC[p,x^\pm]^r$, $\P\subseteq \CC^k$, and $\X\subseteq (\CC^*)^n$, we consider from now on the following statements:

\smallskip
\begin{itemize}[widest=(degAllG),itemindent=*,leftmargin=*]
\item[{\rm (deg1)}]  $F_p$ has a nondegenerate zero  in $(\CC^*)^n$ for some $p\in \CC^k$.
\item[{\rm (degX1)}] $F_p$ has a nondegenerate zero  in $\X$ for some $p\in \P$.
\item[{\rm (degXG)}] There is a nonempty Zariski open subset $\mathcal{U}\subseteq \I$ such that for all $(p,x)\in \mathcal{U}\cap (\P\times \X)$, $x$ is a nondegenerate zero of $F_p$. 
\item[{\rm (degAllG)}] For generic $p\in {\Z}$,  all zeros of $F_p$ in $(\CC^*)^n$ are nondegenerate.
\item[{\rm (setE)}] $\Z_{\P}(\X)$ has nonempty Euclidean interior in $\P$.
\item[{\rm (setZ)}] $\Z_{\P}(\X)$ is Zariski dense in $\CC^k$. 
\item[{\rm {(flatG)}}] {The projection $\pi\from\I\to\CC^k$ is flat at generic $(p,x)\in\I$.}
\item[{\rm (dim1)}] $\VV_{\CC^*}(F_p)$ has pure dimension $n-r$ for at least one $p\in\Z_{\P}(\X)$.
\item[{\rm (dimG)}] $\VV_{\CC^*}(F_p)$ has pure dimension $n-r$ for generic $p\in\Z_{\P}(\X)$.
\item[{\rm (dimX)}] $\VV_{\CC^*}^{\X}(F_p)$ has pure dimension $n-r$ for generic $p\in\Z_{\P}(\X)$.
\item[{\rm (rad)}] $F_p$ generates a radical ideal in $\CC[x^{\pm}]$ for generic $p\in\CC^{k}$, and $F$ generates a radical ideal in $\CC(p)[x^{\pm}]$.
\item[{\rm (reg)}]   $\VV_{\CC^*}(F_p)$ is a nonsingular complex algebraic variety for generic $p\in \CC^k$.
\item[{\rm (real)}] If  $F\in\RR[p,x^\pm]^r$, $\P\subseteq \RR^k$ and $\X\subseteq (\RR^*)^n$, then $\VV_{\RR^*}^{\X}(F_p)$ and $\VV_{\RR^*}(F_p)$ have pure dimension $n-r$ for generic $p\in\Z_{\P}(\X)$.
     \end{itemize}
\smallskip

Recall that $n-r$ is the expected dimension of $\VV_{\CC^*}(F_p)$, while $k+n-r$ is the expected dimension of $\I$. The following main theorem gathers the conclusions of our results so far. 

\begin{theorem}\label{thm:summary_of_implications}
Let $F\in \CC[p,{x}^{\pm}]^r$, $\P\subseteq \CC^k$ locally Zariski dense and $\X\subseteq (\CC^*)^n$.
Assume that the incidence variety $\I$ is irreducible of dimension $k+n-r$ and that 
$\I\cap (\P\times \X)$ is Zariski dense in $\I$. 
The following implications hold:
\begin{enumerate}[label=(\roman*)]
\item  {\rm (deg1)}, {\rm (degX1)} and {\rm (degXG)} are equivalent, and {\rm (degAllG)} implies any of these.
If in addition $r=n$, then these four statements are all equivalent. 
\item   {\rm (setZ)}, {\rm (dimG)} and  {\rm (dim1)} are equivalent.
\item   {\rm (setE)} implies {\rm (setZ)}.
\item   {\rm (deg1)} implies {\rm (setZ)}. 
\item {\rm (dimG)} implies {\rm (dimX)}.
\item  {\rm (degAllG)} implies  {\rm (rad)}, {\rm (reg)} and {\rm (real)}.
\item  {{\rm (flatG)} implies {\rm (setZ)}, and the converse is true if $\I$ is nonsingular. }
\end{enumerate}
\end{theorem}
\begin{proof} Note throughout that by hypothesis, $\Z_\P(\X)\neq \varnothing$.  Statement \textbf{(iii)} is a consequence of \Cref{lem:nonemptyinterior} and statement \textbf{(iv)} follows from \Cref{prop:degdim}(i) and \Cref{lem:denseness_of_ZPX}. 

We show next \textbf{(i)}. We have that {\rm (degXG)}$\Rightarrow${\rm (degX1)}, as $\mathcal{U}$ is a nonempty Zariski open set and $\I\cap (\P\times \X)$ is Zariski dense in $\I$, guaranteeing that the intersection $ \mathcal{U}\cap (\P\times \X)$ is nonempty. The implication {\rm (degX1)}$\Rightarrow${\rm (deg1)} is immediate, {\rm (deg1)}$\Rightarrow${\rm (degXG)} follows from 
\Cref{prop:propoerties_of_nondegenerate_zeros}(iii),  and
{\rm (degAllG)} trivially implies  {\rm (deg1)}. \Cref{thm:nondegF} gives    {\rm (deg1)}$\Rightarrow${\rm (degAllG)} when $r=n$ as by \textbf{(iv)}, 
$\Z$ is Zariski dense in $\CC^n$ and hence contains a nonempty Zariski open set of $\CC^k$.
 
For \textbf{(ii)}, 
for  the implication (setZ)$\Rightarrow$(dimG), \Cref{lem:denseness_of_ZPX} gives that $\Z_{\P}$ is Zariski dense whenever $\Z_{\P}(\X)$ is, and hence (dimG) follows from \Cref{thm:generic_dimension}(ii).
The implication (dimG)$\Rightarrow$(dim1) is immediate. 
Finally, 
if (dim1) holds, then \Cref{thm:generic_dimension} gives that $\Z_{\P}$ is Zariski dense, and by \Cref{lem:denseness_of_ZPX} so is $\Z_{\P}(\X)$, giving (setZ).
 
For \textbf{(v)}, the implication holds as $\VV_{\CC^*}^{\X}(F_p)$ is the nonempty union of irreducible components of $\VV_{\CC^*}(F_p)$ if $p\in \Z_{\P}(\X)$. 

For \textbf{(vi)}, the implications from {\rm (degAllG)} to (real) and (reg) follow from \Cref{prop:propoerties_of_nondegenerate_zeros}. 
 The implication {\rm (degAllG)}$\Rightarrow$(rad)  follows from  \Cref{prop:generic_radicality_in_Laurent_polynomial_ring}. 

{Finally, \textbf{(vii)} is the content of \Cref{prop:flatness}.}
\end{proof}

\newcommand{\implications}{}

 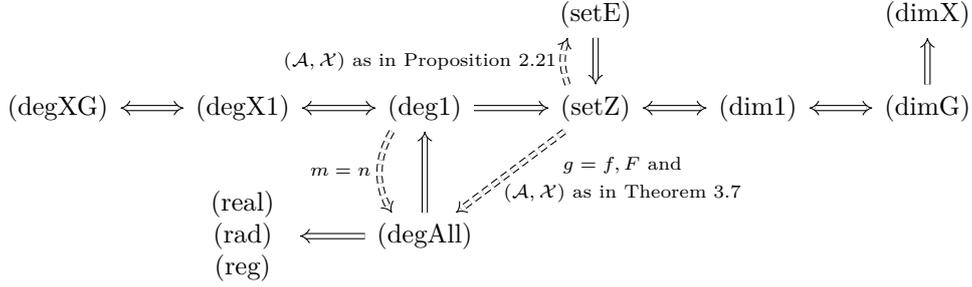
\begin{figure}[t]
    \centering
{\[\begin{tikzcd}
	& &&\mathrm{(setE)} & {\mathrm{(flatG)}} &  {\mathrm{(dimX)}} \\
	 {\mathrm{(degXG)}} & {\mathrm{(degX1)}} & {\mathrm{(deg1)}} 
	& {\mathrm{(setZ)}} & {\mathrm{(dim1)}} & {\mathrm{(dimG)}} \\ 
& {\begin{array}{c}\mathrm{(real)}\\\mathrm{(rad)}\\\mathrm{(reg)}\end{array}} & {\mathrm{(degAllG)}} & & & 
	\arrow[Rightarrow, from=3-3, to=2-3]
	\arrow[Rightarrow, from=3-3, to=3-2]
    \arrow[Rightarrow, to=3-3, from=2-3, dashed, "{\text{$r=n$\hspace{3pt}}\hspace{0.1em}}"', shift right=3,bend right=30]
	\arrow[Leftrightarrow, from=2-1, to=2-2]
	\arrow[Leftrightarrow, from=2-2, to=2-3]
	\arrow[Rightarrow, from=2-3, to=2-4]
		\arrow[Leftrightarrow, from=2-4, to=2-5]
	\arrow[Leftrightarrow, from=2-5, to=2-6]
	\arrow[Rightarrow, from=2-6, to=1-6]
 \arrow[Rightarrow, from=1-4, to=2-4]
	\arrow[Rightarrow, dashed, from=2-4, to=3-3, "\hspace{-10pt}\text{\begin{tabular}{l}
	$F$ augmented vertical \\[-0.3em] and $(\P,\X)$ as in \Cref{thm:main_theorem_for_vertical_systems}\end{tabular}}" very near start]
    \arrow[Rightarrow, to=1-4, from=2-4, dashed, shift left=3, bend left=20, "\text{$(\P,\X)$ as in \Cref{prop:properties_of_alg_defined_pairs}}\hspace{3pt}"]
    \arrow[Rightarrow, from=1-5, to=2-4, shift left=-4pt]
    \arrow[Rightarrow, dashed, from=2-4, to=1-5, "\hspace{2pt}\text{$\I$ nonsingular}"' very near end, bend right=20, start anchor={[yshift=-4pt]north east}] 
\end{tikzcd}\]}
    \caption{Graphical illustration of the implications in \Cref{thm:summary_of_implications}. The dashed arrows indicate implications that hold under additional assumptions. }
    \label{fig:summary}
\end{figure}

The condition that $\I\cap (\P\times \X)$ is Zariski dense in $\I$ only imposes very mild conditions on $\X$, as the following small example illustrates.

\begin{example}\label{ex:restriction_to_rationals} Let $f=x-p$, $\P=\CC$ and $\X=\QQ^*$.
Then $\I=\{(p,p) : p\in \CC^*\}$. The intersection $\I\cap (\P\times \X)=\{(p,p) : p\in \QQ^*\}$ is dense in $\I$ and $\Z_{\P}(\X)=\QQ^*$ is dense in $\CC$.  As (setZ) holds, \Cref{thm:summary_of_implications} applies and the generic dimension is $0$. Note that  (setE) does not hold, showing that the implication (setZ)$\Rightarrow$(setE) is not true for general $F$ and  $\X$. 
\end{example}

\begin{remark}\label{rk:extension_C}
We have chosen to focus on the complex torus because removal of zero coordinates will be necessary in the systems studied in the next section. However, for systems with nonnegative exponents $F\in \CC[p,x]^r$, if the \term{affine incidence variety} 
\[\I_{\CC}:=\{(p,x)\in\CC^{k}\times\CC^n:F_p(x)=0\}\]
is irreducible, \Cref{thm:generic_dimension} holds also in the affine setting over $\CC$ (rather than $\CC^*$) after replacing $\VV_{\CC^*}(F_p)$ by $\VV_{\CC}(F_p)$. 
Additionally, \Cref{lem:denseness_of_ZPX}, \Cref{prop:propoerties_of_nondegenerate_zeros,prop:degdim,prop:generic_radicality_in_Laurent_polynomial_ring}, as well as \Cref{thm:nondegF} extend easily to $\CC$ 
by allowing $\X\subseteq \CC^n$. We conclude that \Cref{thm:summary_of_implications} holds in the affine setting over $\CC$ as long as $\I_{\CC}$ is irreducible and $\I_{\CC}\cap (\P\times \X)$ is Zariski dense in $\I_{\CC}$. 
\end{remark}

In \Cref{sec:vertical_systems}, we will see that for augmented  vertically parametrized systems, the {nine} first statements (deg1)--(dimG) from the beginning of the subsection are equivalent.
Thus, all the desired properties of the polynomial system will rely on the existence of  a nondegenerate zero in the complex torus.
The key part of the argument will be to show that
$\Z$ has nonempty Euclidean interior if and only if $F$ has a nondegenerate zero. This phenomenon {need not} happen for general families, as example \eqref{eq:example_nonlinear_coefficients} shows. 

\subsection{Zariski denseness and nonempty Euclidean interior}
For the stronger result discussed at the end of the previous subsection to hold, we will need to restrict to sets $\P$ and $\X$ where the converse of \Cref{lem:nonemptyinterior} holds, that is (setZ) implies (setE) (c.f. \Cref{ex:restriction_to_rationals}). In this final subsection, we identify pairs $(\P,\X)$ for which this holds.

\noindent
\begin{minipage}{\textwidth}
\begin{lemma}\label{lemma:denseness_of_constructible_semialgebraic}\leavevmode
\begin{enumerate}[label=(\roman*)]
    \item A constructible set $S\subseteq\CC^k$ is Zariski dense in $\CC^k$ if and only if it has nonempty Euclidean interior in $\CC^k$.
    \item A semialgebraic set $S\subseteq\RR^k$ is Zariski dense in $\CC^{k}$ if and only if it has nonempty Euclidean interior in $\RR^k$.
\end{enumerate}
\end{lemma}
\end{minipage}

\begin{proof}
For part (i),  as  $S$ is constructible,   $S=\bigcup_{i=1}^t Z_i\cap U_i$, for some irreducible Zariski closed sets $Z_i\subseteq\CC^k$ and nonempty Zariski open subsets $U_i\subseteq\CC^k$. The statement now follows from the fact that the Zariski closure satisfies $\overline{S}=\bigcup_{i=1}^t Z_i$. 
In case (ii),  \cite[Proposition~2.8.2]{BCR} gives that $S$ is Zariski dense in $\CC^k$ if and only if its semialgebraic dimension is $k$, which in turn is equivalent to $S$ having nonempty Euclidean interior in $\RR^k$.
\end{proof}

In order to apply \Cref{lemma:denseness_of_constructible_semialgebraic}, we note that $\Z_{\P}(\X)$ is constructible or semialgebraic, if both $\P$ and $\X$ are defined by algebraic data in the sense of the following definition.

\begin{definition}\label{def:algebraically_defined_pair}
A pair of sets $(\P,\X)$ where $\P\subseteq\CC^k$ and $\X\subseteq \CC^n$ is
said to be \term{algebraically defined} if  
it satisfies one of the following conditions:
\begin{enumerate}[label=(\roman*)]
\item $\P$ is a Zariski open subset of $\CC^k$ and $\X$ is constructible.
\item $\P\subseteq\RR^k$ and $\X\subseteq \RR^n$ are both semialgebraic, and $\P$ is locally Zariski dense.
\item $\P\subseteq\RR^{k}$ is semialgebraic and  locally Zariski dense, and $\X\subseteq \CC^n$ is constructible.  
\end{enumerate}
\end{definition}

Examples of algebraically defined pairs include $(\RR^k,\RR^n_{>0})$, $(\RR^k,(\RR^*)^n)$, and $(\CC^k,(\CC^*)^n)$.

\begin{proposition}\label{prop:properties_of_alg_defined_pairs}
Let $F\in \CC[p,x^{\pm}]^r$ and   $(\P,\X)\subseteq \CC^k \times (\CC^*)^n$ be an algebraically defined pair. Then:
 \begin{enumerate}[label=(\roman*)]
\item $\Z_{\P}(\X)$ is either a semialgebraic subset of $\RR^k$ or a constructible subset of $\CC^k$.
\item (setZ) and (setE) in \Cref{thm:summary_of_implications} are equivalent.
\end{enumerate}
\end{proposition}

\begin{proof}
For (i),   in
case \Cref{def:algebraically_defined_pair}(i), the set $\Z_{\P}(\X)\subseteq \CC^k$ is constructible by Chevalley's theorem.  
For the case \Cref{def:algebraically_defined_pair}(ii), the Tarski--Seidenberg Theorem \cite[Theorem~2.2.1]{BCR} gives that $\Z_{\P}(\X)\subseteq\RR^k$ is semialgebraic.
Finally, in case \Cref{def:algebraically_defined_pair}(iii), Chevalley's theorem gives that $\Z_{\CC^{k}}(\X)\subseteq\CC^{k}$ is constructible. This in turn implies that \[\Z_{\P}(\X)=\Z_{\CC^{k}}(\X)\cap\P=(\Z_{\CC^{k}}(\X)\cap\RR^{k})\cap\P\] is an intersection of semialgebraic sets and therefore semialgebraic (note that the real points of a constructible set in $\CC^{k}$ form a semialgebraic set in $\RR^{k}$). 

For (ii), the reverse implication is \Cref{lem:nonemptyinterior} as $\P$ is locally Zariski dense. 
For the forward implication, by \Cref{lemma:denseness_of_constructible_semialgebraic}, if $\Z_{\P}(\X)$ is Zariski dense in $\CC^k$, then it has nonempty Euclidean interior in $\RR^k$ in cases (ii) and (iii) of \Cref{def:algebraically_defined_pair} and in $\CC^k$ in case (i). This in turn implies that $\Z_{\P}(\X)$ has nonempty Euclidean interior in $\P$, concluding the proof. 
\end{proof}

\begin{remark}
The condition of $\P$ being locally Zariski dense cannot be replaced by $\P$ being Zariski dense in \Cref{lem:nonemptyinterior} and in \Cref{def:algebraically_defined_pair}(ii),(iii). 
The fact that $\Z_{\P}(\X)$ is  semialgebraic, gives that $\Z_{\P}(\X)$ being Zariski dense in $\CC^k$ is equivalent to having nonempty Euclidean interior in $\RR^k$. The latter might not be equivalent to having nonempty Euclidean interior in $\P$. For example, this fails for the semialgebraic set $\P=\RR^2_{>0} \sqcup \{p_1+p_2=0\}$ if $\Z_{\P}(\X)$ is contained in the line with equation  $p_1+p_2=0$. This phenomenon cannot happen if, with the Euclidean topology, the closure of $\P$ equals the closure of the  interior of $\P$ in $\RR^k$, as $\P$ is then locally Zariski dense. The latter will typically be the case in applications. 
\end{remark}

\section{Vertically parametrized systems} 
\label{sec:vertical_systems}

We now turn our attention to vertically parametrized systems. We study their incidence varieties and the set of parameter values for which all zeros are degenerate, and use this to prove our main result \Cref{thm:main_theorem_for_vertical_systems}. We end the section by commenting on how our results specialize for freely parametrized systems, as well as on conditions that ensure that our results extend from $\CC^*$ to $\CC$.

\subsection{Structure of the systems}
By a \term{vertically parametrized system} (or \term{vertical system} for short), we mean a parametric system
$F\in\CC[a,x^\pm]^s$ for $s\leq n$, 
with parameters $a=(a_1,\ldots,a_m)$ and variables $x=(x_1,\ldots,x_n)$, with the following properties:
\begin{itemize}
    \item The coefficients of the monomials in $x$ are homogeneous linear forms in $a$.
    \item Each $a_i$ appears in at most one column of the coefficient matrix of $F$ regarded in $\CC(a)[x^\pm]^s$.
    \item The coefficient matrix of $F$ regarded in $\CC[a,x^\pm]^s$ has full rank. 
\end{itemize}
This is equivalent to the existence of matrices 
$C\in \CC^{s\times m}$ and $M\in \ZZ^{n\times m}$ such that
\begin{equation}\label{eq:vertical_system}
F = C (a \star x^M), \qquad \rk(C)=s\leq n.
\end{equation}
The rows of $C$ produce $s$ linear combinations of $m$ monomials encoded by the columns of $M$, where the $i$-th monomial is scaled by the parameter $a_i$. Note that $M$ might have repeated columns; this corresponds to the same monomial appearing several times in the system, accompanied with different parameters. 

More generally, we will consider \term{augmented vertically parametrized systems} of the form
\[\left(C(a\star x^M),\, Lx-b\right)\in\CC[a,b,x^\pm]^{s+\ell}, \qquad s+\ell \leq n,\]
where we also include $\ell$ affine equations, encoded by a coefficient matrix $L\in\CC^{\ell\times n}$ 
and parametric constant terms $b=(b_1,\ldots,b_\ell)$. Geometrically, this corresponds to intersecting the variety given by $C (a \star x^M)$ by a   parallel translate of $\ker(L)$.
We observe that as $C$ has full row rank, the coefficient matrix of the system $F$, regarded in $\CC[a,b,x^\pm]^s$, has full rank, independently of the rank of $L$ (c.f. \Cref{rem:dim}). 
When $\ell=0$, an augmented vertical system is simply a vertical system. 

Similarly to what we have done in the previous section, we will restrict the parameter values to some sets $\A\subseteq\CC^m$ and $\B\subseteq\CC^\ell$, and the variable values to some set $\X\subseteq(\CC^*)^n$. In the notation of the previous section, the parameter space, number of parameters, and number of polynomials in the system are, respectively,
\begin{equation}\label{eq:newnotation}
\P=\A\times\B \subseteq \CC^{m} \times \CC^\ell, \qquad k=m+\ell \qquad \text{and}\qquad r=s+\ell. 
\end{equation}

\subsection{The incidence variety}
We now show that the incidence varieties for augmented vertical systems  are irreducible  and nonsingular everywhere, and derive some basic facts about their geometry. In particular, it follows that \Cref{thm:generic_dimension} is applicable.

To this end, with $h^{-1}=(h_1^{-1},\ldots,h_n^{-1})$, we consider the map
\begin{align}\label{eq:parametrization_of_incidence_variety}
\phi \colon \ker(C)\times (\CC^*)^n \rightarrow \CC^{m}\times\CC^{\ell}\times (\CC^*)^n,\quad (w,h)\mapsto (w\star h^M,Lh^{-1},h^{-1}). 
\end{align}

\begin{theorem}\label{thm:irreducible_smooth_incidence_variety}
Let $F=(C(a\star x^M),Lx-b)$ be an augmented vertical system with $C\in \CC^{s\times m}$ of full rank $s$, $L\in \CC^{\ell\times n}$  with $s+\ell\leq n$ and $M\in \ZZ^{n\times m}$. Then the following statements hold:
 \begin{enumerate}[label=(\roman*)]
\item  The map $\phi$ is injective and $\I = \im(\phi)$.
\item $\I$  is irreducible  of dimension $m+n-s$ and has no singular points.
\end{enumerate}
\end{theorem}

\begin{proof}
To prove (i),  injectivity is straightforward, and for surjectivity onto $\I$, we first note that if $(w,h)\in \ker(C) \times (\CC^*)^n$, then 
\[F(w\star h^M,Lh^{-1},h^{-1})=(C((w\star h^M)\star (h^{-1})^M),Lh^{-1}-Lh^{-1})=(C(w),0)=(0,0).\]
Hence, $\im(\phi)\subseteq\I$. To show the reverse inclusion, given $(a,b,x)\in\I$, we have that  $a\star x^M\in \ker(C)$.  
By letting $w=a\star x^M$ and $h=x^{-1}$, we obtain $\phi(w,h)=(a,b,x)$.

For part (ii), the irreducibility and dimension claims follow from the isomorphism of varieties $\I\cong\ker(C)\times (\CC^*)^n$ from  (i). To prove nonsingularity, we observe that each $(a,b,x)\in\I$ is a nondegenerate zero of $F$. Indeed, the Jacobian of $F$ has the form
\[J_F(a,b,x)=\begin{bmatrix} C \diag(x^M) & 0_{s\times\ell} & * \\0_{\ell\times m} & -{\rm Id}_\ell & L  \end{bmatrix}.\]
As both $C \diag(x^M)$ and ${\rm Id}_\ell$ have maximal row rank for $x\in(\CC^*)^n$, so has $J_F(a,b,x)$.
\end{proof}

\subsection{Nondegenerate zeros and generic consistency}
It follows from \Cref{prop:degdim} that if an augmented vertical system $F=(C(a\star x^M),Lx-b)$ has a nondegenerate zero, then $\Z$ has nonempty Euclidean interior in the parameter space $\CC^{m}\times\CC^{\ell}$. 
We will see next that the converse also holds, that is, if $\Z$ has nonempty Euclidean interior, then necessarily 
the system $F$ has a nondegenerate zero for some $(a,b)\in\CC^m\times\CC^\ell$. 
We start by deriving a parametrization of $\Z$ from \eqref{eq:parametrization_of_incidence_variety}, and proceed to study how degenerate zeros behave with respect to this parametrization. This, together with Sard's lemma, gives the main result \Cref{thm:degenerate_zeros_of_vertical_systems} of this subsection.

Let $G\in\CC^{m\times (m-s)}$ be  Gale dual to $C$, in the sense that the columns of $G$ form a basis for $\ker(C)$.
The parametrization $\phi$ in \eqref{eq:parametrization_of_incidence_variety}  gives rise to a parametrization of $\Z$ via the Laurent polynomial map
\[\Phi:=\pi\circ\phi\circ(G\times\id)\from\CC^{m-s}\times(\CC^*)^n\to\CC^m\times\CC^\ell, \qquad (u,h)\mapsto ((Gu)\star h^M,Lh^{-1}),\]
where $\id$ is the identity on $(\CC^*)^n$, and $\pi\from\CC^{m}\times\CC^\ell\times (\CC^*)^n\to\CC^m\times\CC^\ell$ the canonical projection.

Consider now the set of parameters that give rise to degenerate zeros:
\begin{equation}\label{eq:D}
    \D := \{(a,b)\in\CC^m\times\CC^\ell : \text{$F_{a,b}$ has a degenerate zero in $(\CC^*)^n$}\}\subseteq \Z.
\end{equation}

We introduce the following matrix for $(w,h)\in\CC^{m-s}\times(\CC^*)^n$
\begin{align}\label{eq:Q}
Q(w,h)&:=\begin{bmatrix}C \diag(w)M^\top \diag(h) \\ L\end{bmatrix}\in \CC^{(s+\ell)\times n},
\end{align}
and the set
\[\Delta := \{(u,h)\in \CC^{m-s}\times  (\CC^*)^n :\rk(Q(Gu,h))<s+\ell\}.\]

\begin{proposition}\label{prop:description_of_Z_and_D}
Let $F=(C(a\star x^M),Lx-b)$ be an augmented vertical system and consider the notation above. Then the following statements hold:
\begin{itemize}
\item[(i)] Let $w\in \ker(C)$, $x,h\in (\CC^*)^n$, $a\in \CC^m$ and $b\in \CC^{\ell}$ satisfy 
$(a,b,x)=\phi(w,h)$. Then   
\[J_{F_{a,b}}(x)=Q(w,h). \]
\item[(ii)] $F_{a,b}$ has a nondegenerate zero in $(\CC^*)^n$  for some $(a,b)\in \CC^m\times\CC^\ell$ if and only if $\rk(Q(w,h))=s+\ell$ for some $w\in \ker(C)$ and $h\in (\CC^*)^n$.
\item[(iii)] $\Z = \im(\Phi)$ and $\D = \Phi(\Delta)$. 
\end{itemize}
\end{proposition}
\begin{proof}
An easy computation shows that the Jacobian of $F_{a,b}$ is given by
\[J_{F_{a,b}} ({x}) =\begin{bmatrix} C \diag({a}\star {x}^M)M^\top\diag(x^{-1})  \\ L \end{bmatrix}, \]
from which (i) follows. Statement (iii) follows from part (i) and \Cref{thm:irreducible_smooth_incidence_variety}(i).
The   description of $\Z$ in (ii) follows  from \Cref{thm:irreducible_smooth_incidence_variety}(i), and that of $\D$  follows from this and part (i). 
\end{proof}

\begin{proposition}\label{lem:critical_points}
With the notation above,  $\Delta$ agrees with the set of critical points of $\Phi$ (equivalently, $J_\Phi(u,h)$ does not have full rank if and only if $(u,h)\in \Delta$).
\end{proposition}

\begin{proof}
We start by noting that $(u,h)\in\CC^{m-s} \times (\CC^*)^n$ is a critical point of $\Phi$ if and only if 
$\dim(\ker(J_{\Phi}(u,h)))>n-(s+\ell)$, and that  $(u,h)\in \Delta$ if and only if $\dim(\ker(Q(Gu,h)))>n-(s+\ell)$. Hence, all we need is to show that 
$\dim(\ker(J_{\Phi}(u,h))) = \dim(\ker(Q(Gu,h)))$. 

A simple computation gives that
\[
J_{\Phi}(u,h)=\begin{bmatrix}
\diag(h^M) G & \diag(Gu \star {h}^M)M^\top\diag(h^{-1})\\
0 & -L\diag(h^{-2})\end{bmatrix}, 
\]
and hence its rank agrees with that of
\[P(u,h)=\begin{bmatrix}
G & \diag(Gu)M^\top\diag(h)\\
0 & L\end{bmatrix} \in \CC^{(m+\ell) \times (m+n-s)}.\] 
As $CG=0$ we further have that 
\begin{equation}\label{eq:RQ} 
\begin{bmatrix}C & 0\\ 0 & {\rm Id}_\ell\end{bmatrix}P(u,h) = \begin{bmatrix}
     0_{(s+\ell) \times (m-s)} &  Q(Gu,h) \end{bmatrix}
\end{equation}
and
\[\ker\Big(\begin{bmatrix}C & 0\\ 0 & {\rm Id}_\ell\end{bmatrix}\Big)\cap\im(P(u,h)) = \im(G) \times \{0\}. \] 
From this and \eqref{eq:RQ} we obtain
\begin{align*}
  \dim (\ker (P(u,h)) &= \dim\Big(\ker\Big(\begin{bmatrix}C & 0\\ 0 & {\rm Id}_\ell\end{bmatrix}P(u,h)\Big)\Big) - \dim \Big(\ker\Big(\begin{bmatrix}C & 0\\ 0 & {\rm Id}_\ell \end{bmatrix}\Big)\cap\im(P(u,h))\Big) \\
  &=  \dim (\ker (Q(Gu,h)) + (m-s) - (m-s) =  \dim (\ker (Q(Gu,h)).\qedhere
\end{align*}
\end{proof}

The next statement is given for zeros in $(\CC^*)^n$ and for the full parameter space  $\CC^{m}\times\CC^\ell$. We will later see in \Cref{thm:main_theorem_for_vertical_systems} that the statement also holds for more general sets $\A$, $\B$ and $\X$.

\begin{proposition}\label{prop:D_contained_in_hypersurface}
The set $\D$ in \eqref{eq:D} is contained in a hypersurface of $\CC^{m}\times\CC^\ell$.  
\end{proposition}
\begin{proof}
As $\D=\Phi(\Delta)$ by \Cref{prop:description_of_Z_and_D}(iii), it follows from \Cref{lem:critical_points} that $\D$ is the set of critical values of $\Phi$. Now Sard's Lemma (see, e.g., \cite[Proposition~B.2]{S-E-D_S} and \cite[Proposition~3.7]{Mumford}), adapted to very affine varieties, gives that $\dim(\overline{\D})<m+\ell$.
\end{proof}

We summarize our conclusions on degenerate zeros of augmented vertical systems as follows.
 
\begin{theorem}\label{thm:degenerate_zeros_of_vertical_systems}
Let $F=(C(a\star x^M),Lx-b)$ be an augmented vertical system  with $C\in \CC^{s\times m}$ of full rank $s$, $L\in \CC^{\ell\times n}$ with $s+\ell\leq n$ and $M\in \ZZ^{n\times m}$.  The following statements are equivalent:
\begin{enumerate}[label=(\roman*)]
\item For all $(a,b)\in \Z$, all zeros of  $F_{a,b}$ in $(\CC^*)^n$ are degenerate.
\item For every $(a,b)\in \Z$, it holds that $F_{a,b}$ has a degenerate zero in $(\CC^*)^n$, i.e., $\D=\Z$.
 \item $\Z$ has empty Euclidean interior in $\CC^m\times\CC^\ell$.
\end{enumerate}
Additionally, if $\Z$ has nonempty Euclidean interior in $\CC^m\times\CC^\ell$, then all zeros of $F_{a,b}$ are nondegenerate for generic choices of $(a,b)\in \CC^m\times\CC^\ell$.
\end{theorem}

\begin{proof}  
The implication (i)$\Rightarrow$(ii) holds trivially, while (ii)$\Rightarrow$(iii) is a consequence of \Cref{prop:D_contained_in_hypersurface}. The implication
(iii)$\Rightarrow$(i) follows by \Cref{prop:degdim}. For the last statement, it is enough to consider the Zariski open subset $\mathcal{U}=\Z \setminus \overline{\D}$ of $\Z$, which is nonempty by \Cref{prop:D_contained_in_hypersurface}.
\end{proof}

\begin{remark}
For $\mathbb{G}=\RR^*$ and $\mathbb{G}=\RR_{>0}$,  
it is straightforward to verify that restricting the map
$\phi\from\ker(C)\times(\CC^*)^n\to\I$ to $(\ker(C)\cap\mathbb{G}^{m})\times\mathbb{G}^n$ provides a parametrization of $\I \cap ((\mathbb{G}^{m}\times \RR^\ell) \times \mathbb{G}^n)$, which after projection gives a parametrization of 
$\Z_{\mathbb{G}^{m}\times \RR^\ell}(\mathbb{G}^n)$. 
\end{remark}

\subsection{The main theorem on generic dimension}

We now apply the results of the previous subsections to complete the equivalences in \Cref{thm:summary_of_implications}, and thereby unify \Cref{thmalph:main_result_over_C} and \Cref{thmalph:main_result_over_R} from the introduction.
Recall from \eqref{eq:newnotation} the relation between the notation in \Cref{thm:summary_of_implications} and the current section, and recall the matrix $Q(w,h)$ defined in \eqref{eq:Q}. 
 
\begin{theorem}\label{thm:main_theorem_for_vertical_systems}
Let $F=(C(a\star x^M),Lx-b)$ be an augmented vertical system with $C\in \CC^{s\times m}$ of full rank $s$, $L\in \CC^{\ell\times n}$ with $s+\ell\leq n$ and $M\in \ZZ^{n\times m}$. Let $\P=\A\times\B\subseteq \CC^m \times \CC^\ell$ and suppose that $(\P,\X)$ is an algebraically defined pair of sets, for which  
$\I\cap (\P\times \X)$ is Zariski dense in $\I$. 
Then the following holds:
\begin{enumerate}
\item[(i)] The statements {\rm (deg1)}, {\rm (degX1)}, {\rm (degXG)}, {\rm (degAllG)}, {\rm (setE)}, {\rm (setZ)}, {{\rm (flatG)}},  {\rm (dim1)} and  {\rm (dimG)} are all equivalent to the condition 
\[\rk(Q(w,h))=s+\ell\quad\text{for some $(w,h)\in \ker(C) \times (\CC^*)^n$}.\]
\item[(ii)] Any of the statements mentioned above implies {\rm (dimX)}, {\rm (reg)} and {\rm (real)}. 
\item[(iii)]  The statement {\rm (rad)} holds, independently of the other statements. 
\end{enumerate}
\end{theorem}
\begin{proof}
The equivalence between (deg1) and the matrix rank condition is \Cref{prop:description_of_Z_and_D}(ii). {As $\I$ is irreducible and nonsingular of dimension $m+n-s$, the assumptions for \Cref{thm:summary_of_implications} are satisfied, and there is an equivalence between (setZ) and (flatG)}. Therefore, for (i) and (ii), it suffices to show that (setZ) implies (degAllG) and (setE). The equivalence between (setZ) and (setE) is  \Cref{prop:properties_of_alg_defined_pairs}(ii). 
By \Cref{lem:denseness_of_ZPX} and \Cref{prop:properties_of_alg_defined_pairs}, $\Z$ has nonempty Euclidean interior in $\CC^m\times\CC^\ell$ if (setZ) holds.  
By \Cref{thm:degenerate_zeros_of_vertical_systems}, this implies (degAllG). Finally, to prove (iii), we note that if (degAllG) holds, so does (rad) by \Cref{thm:summary_of_implications}. If, on the other hand, (degAllG) does not hold, then (i) gives that (setZ) does not hold either, and hence it follows by \Cref{thm:generic_dimension} that $\VV_{\CC^*}(F_{a,b})=\varnothing$ for generic $(a,b)\in \CC^m\times\CC^\ell$. Thus, generically, $\langle F_{a,b}\rangle=\CC[x^{\pm}]$, which is radical.
\end{proof}

\begin{remark}\label{rem:remarks_on_main_result}\leavevmode
\begin{enumerate}[label=(\roman*)]
\item As indicated in \Cref{thm:summary_of_implications}(i),  the implication (deg1)$\Rightarrow$(degAllG) already followed from \Cref{thm:nondegF} for square augmented vertical systems. The results of this section have shown that this   also holds  in the non-square case.
\item For any $F$, $\P$ and $\X$  satisfying the assumptions of \Cref{thm:main_theorem_for_vertical_systems}, any statement in \Cref{thm:main_theorem_for_vertical_systems}(i) is equivalent to (deg1), and hence to any of the same statements for $\P=\CC^m\times\CC^\ell$ and $\X=(\CC^*)^n$ instead.
\item If (degAllG) holds under the assumptions of \Cref{thm:main_theorem_for_vertical_systems}, then we also have that all zeros of $F_{a,b}$ are nondegenerate for generic $(a,b)\in \Z_\P(\X)$.
\item When $\ell=0$, the rank condition in \Cref{thm:main_theorem_for_vertical_systems}(i) simplifies to
\[\rk(C\diag(w)M^\top ) = s \quad \text{for some }w\in \ker(C). \]
\item If a vertical system $C(a\star x^M)$ is generically consistent, then so is the augmented vertical system $(C(a\star x^M),Lx-b)$ for  generic $L\in\CC^{\ell\times n}$ for $\ell\leq n-s$.
\end{enumerate}
\end{remark}

\Cref{thm:main_theorem_for_vertical_systems} tells us that for augmented vertical systems, only one of two extreme scenarios occurs, provided the parameter sets $\A$ and $\B$ and the domain $\X$ of the variables satisfy some mild conditions: Either the system is generically consistent, has generically the expected dimension $n-(s+\ell)$, and generically all zeros are nondegenerate, or the system is generically inconsistent, has never the expected dimension {over $\CC$}, and all zeros are degenerate. 

In \Cref{subsec:computational} we will discuss, and show several examples of, how to computationally check which of these two scenarios we are in, for a given augmented vertical system.

\begin{remark}
The conclusions of \Cref{thm:main_theorem_for_vertical_systems} also hold for parametric systems $F_{\cdot,b}$ given by the restriction of an augmented vertical system $F$ to a fixed value of $b\in L((\CC^*)^n)$.
In this case, the incidence variety is irreducible {and nonsingular} of dimension $n+m-(s+\ell)$, with parametrization
\[\ker(C) \times U_b   \rightarrow \CC^m\times (\CC^*)^n,  \qquad (w,h)  \mapsto (w\star h^M,h^{-1}), \]
where $U_b:=\{ h\in (\mathbb{C^*})^n: Lh^{-1}=b\}$. Using this, the proof of   \Cref{lem:critical_points} (and subsequently \Cref{thm:main_theorem_for_vertical_systems}) can be adapted to this system.
\end{remark}

\begin{remark}[Non-vertical systems]
One of the restrictions we imposed on vertical systems is that each parameter  always multiplies the same monomial. Relaxing this restriction might make \Cref{thm:main_theorem_for_vertical_systems} fail. 
In particular, this more general class of systems might not give rise to irreducible incidence varieties, and hence the implications in \Cref{thm:summary_of_implications} do not necessarily hold. For example, (setZ) and (deg1) might hold, but not (dimG). Additionally, the implication (setZ)$\Rightarrow$(deg1), which holds for vertical systems by \Cref{thm:main_theorem_for_vertical_systems}, might not hold.
This is illustrated by the examples we saw in \eqref{eq:non_vertical_systems1} and \eqref{eq:non_vertical_systems2}.
\end{remark}

\subsection{Computational considerations}\label{subsec:computational}

Consider an augmented vertical system with the following scenarios regarding the sets $\A$, $\B$, and $\X$, which are common in applications:
\begin{enumerate}
\item[(a)] $\B = \RR^{\ell}$ and\\[0.2em] 
$(\A,\X)\in\{ (\RR^m_{>0},\RR^n_{>0}),((\RR^*)^m,\RR^n_{>0}),
((\RR^*)^m,(\RR^*)^n), (\RR^m,\RR^n_{>0}),
(\RR^m,(\RR^*)^n)\}.$
\vspace{0.2em}
\item[(b)] $\A\in\{(\CC^*)^m,\CC^m\}$, $\B=\CC^{\ell}$, and $\X=(\CC^*)^n$. 
\end{enumerate}
For all of these, $\P=\A\times\B$ is locally Zariski dense  and the pair $(\P,\X)$ is algebraically defined. 
Assuming $F$ has real coefficients and $\X\subseteq (\RR^*)^n$, it follows from
\Cref{rk:PxX_dense_in_I} that
$\I\cap (\P\times \X)$ is Zariski dense in $\I$ whenever 
$\I\cap (\P\times \X) \neq \varnothing$. This, in turn, is characterized as follows. 

\begin{proposition}\label{prop:nonempty_intersection_with_incidence_variety}
Let $\A,\B,\X$ be as in (a) or (b) above, and $F$ be an augmented vertical system (with real coefficients in case (a)). Then $\I\cap (\A\times \B \times \X)\neq \varnothing$ if and only if $\ker(C)\cap\A\neq\varnothing$.
\end{proposition}

\begin{proof}
The conditions on $\A$ and $\X$ imply that for all $a\in \A$ and $x\in \X$, it holds $a\star x^M\in \A$. Hence, if $\ker(C)\cap\A=\varnothing$, 
then $\VV_{\CC^*}(F_{a,b})=\varnothing$ for all $(a,b)\in\A\times\B$. On the other hand, if $\ker(C)\cap\A\neq\varnothing$, then $(1,\ldots,1)\in\VV_{\CC^*}(F_{a,b})\cap \X$ for all $a\in \ker(C)\cap\A$, with $b:= L (1,\dots,1)^\top\in \B$. 
\end{proof}

In the special case $\A=\RR_{>0}^m$, checking $\ker(C)\cap\RR_{>0}^m\neq \varnothing$ corresponds to showing the existence of an interior point of the polyhedral cone $\ker(C)\cap\RR_{\geq 0}^m$. If $C$ has rational entries, this is a straightforward computation using linear programming.

The rank condition from \Cref{thm:main_theorem_for_vertical_systems} can be checked in the following way: Let $G\in\CC^{m\times(m-s)}$ be a Gale dual to $C$, whose columns form a basis for $\ker(C)$. We want to check whether there exists some $(u,h)\in\CC^{m-s}\times(\CC^*)^n$ such that
\begin{equation}\label{eq:rank_condition_computational}
\rk\begin{bmatrix}C\diag(Gu)M^\top\diag(h)\\L\end{bmatrix}=s+\ell .
\end{equation}
Equality \eqref{eq:rank_condition_computational} holds for all $(u,h)$ in a Zariski open subset of $\CC^{m-s}\times(\CC^*)^n$, so if this set is nonempty, then \eqref{eq:rank_condition_computational} holds for a randomly chosen $(u,h)$  with probability $1$, given an appropriate probability measure on $\CC^{m-s}\times(\CC^*)^n$. Hence, we pick a random pair $(u,h)$ and compute the rank in \eqref{eq:rank_condition_computational} with exact arithmetic. If the rank is $s+\ell$, then we are in the generically consistent scenario. If not, we can suspect that we are in the generically inconsistent scenario, and to conclusively prove this, we view the matrix in \eqref{eq:rank_condition_computational} as a symbolic matrix with indeterminates $(u,h)$ and verify that all $(s+\ell)$-minors are the zero polynomial. 

In the special case when $\ell=0$, it follows from \Cref{rem:remarks_on_main_result}(iv) that it suffices to check whether there is some $u\in\CC^{m-s}$ such that
\begin{equation}\label{eq:rank_condition_computational_vertical}
    \rk(C\diag(Gu)M^\top)=s.
\end{equation}

\begin{example}\label{ex:examples_from_intro_revisited}
Consider the vertical system \eqref{eq:vertical_example} from the introduction. Using the matrices
\[C=\begin{bmatrix}1&3&1&0
\\0&1&2&1
\end{bmatrix},\quad 
M=\begin{bmatrix}2&2&0&1\\ 2&0&0&0\\ 0&1&1&1\\ 0&0&2&1\end{bmatrix}, \quad 
G=\begin{bmatrix}
5&3\\ -2&-1\\ 1&0\\ 0&1
\end{bmatrix},\]
we obtain
\[ C\diag(Gu)M^\top =\begin{bmatrix}
-2 u_{1} & 10 u_{1}+6 u_{2} & -5 u_{1}-3 u_{2} & 2 u_{1} 
\\
 -4 u_{1}-u_{2} & 0 & 0 & 4 u_{1}+u_{2} 
\end{bmatrix},
  \]
  and we see that \eqref{eq:rank_condition_computational_vertical} holds for $u=(1,1)$. Hence, the system \eqref{eq:vertical_example} is generically consistent. Consider now the augmented vertical system \eqref{eq:augmented_vertical_example}, with $C$ and $M$ as above, and 
\[L=\begin{bmatrix}
0 & 1 & 2 & 0\\ 0 & 1 & 1 & 0\\
\end{bmatrix}.\]
In this case, the matrix in condition \eqref{eq:rank_condition_computational} is square, and its determinant is the zero polynomial, when $(u,h)$ are viewed as variables. We conclude that  system \eqref{eq:augmented_vertical_example} is generically inconsistent. 
(We do, however, obtain a generically consistent system for any choice of $L\in\CC^{2\times 4}$ with full rank, such that column 1 or 4, as well as column 2 or 3, have at least one nonzero entry.)
\end{example}

\begin{example}\label{ex:critical_points_revisited}
Going back to the application to critical points in optimization mentioned in the introduction, we consider bivariate polynomials of the form $f=a_1 x_1 + a_2 x_1x_2 + a_3 x_2^2$. The matrix $M$ in the square vertical system in \eqref{eq:critical_system} and a Gale dual matrix are
\[M = \begin{bmatrix}
    1 & 1 & 0 \\ 0 & 1 & 2
\end{bmatrix}, \qquad G = \begin{bmatrix}
    2  \\ -2 \\ 1
\end{bmatrix}. \]
This gives that
\[M\diag(Gu )M^\top = \begin{bmatrix}
    0 & -2u \\-2u & 2u
\end{bmatrix}, \]
which has rank $2$ for all $u\neq 0$. 
As the rank condition in \eqref{eq:rank_condition_critical_points} holds, $f_a$ has a finite, positive number of critical points in $(\CC^*)^2$ for generic $a\in\CC^3$. If we now impose the coefficients to be real such that $a_1,a_3>0$ and $a_2<0$, the polynomial $f$
 can be written instead as $f=a_1 x_1 - a_2 x_1x_2  + a_3 x_2^2$ with $a_i>0$. 
The coefficient matrix of the vertical system  encoding the critical points becomes 
 \[C = \begin{bmatrix}
    1 & -1 & 0 \\ 0 & -1 & 2
\end{bmatrix}, \]
which satisfies $\ker(C)\cap\RR_{>0}^3\neq\varnothing$. From this, \Cref{thm:main_theorem_for_vertical_systems} with $\A=\RR^3_{>0}$ and $\X=\RR^2_{>0}$ tells us that $f_a$ has positive critical points for $a$ 
in a subset of $\RR^3_{>0}$ with nonempty Euclidean interior. 
\end{example}

\subsection{Freely parametrized systems}\label{subsec:freely_parametrized}
We now turn our attention to the special case of freely parametrized systems, and record some corollaries of our main results.

Given finite sets $\S_1,\dots,\S_s\subseteq \ZZ^{n}$, we consider the corresponding freely parametrized family
\[\F_\mathrm{free}(\S_1,\dots,\S_s):=\{F=(f_1,\ldots,f_s)\in\CC[x^{\pm}]^s:\text{$\operatorname{supp}(f_i)\subseteq \S_i$ for all $i=1,\ldots,s$}\}\]
and say that a property holds generically in $\F_\mathrm{free}(\S_1,\dots,\S_s)$ if it does so under the isomorphism $\F_\mathrm{free}(\S_1,\dots,\S_s)\cong\prod_{i=1}^s\CC^{\S_i}$ that identifies each polynomial with its coefficients.

The set $\F_\mathrm{free}(\S_1,\dots,\S_s)$ can be encoded as a vertical system of the form \eqref{eq:vertical_system} in the following way. Let $m_i$ be the cardinality of $\S_i$, set $m=m_1+\dots+m_s$, and  define $C\in \CC^{s\times m}$ to be the block diagonal matrix 
\begin{equation}\label{eq:Cfree}
    C = \begin{bmatrix}
    C_1 & 0 & \dots & 0 \\ 0 & C_2 & \dots & 0 \\ \vdots & \vdots & \ddots & 0 \\ 0 & 0 & \dots & C_s
\end{bmatrix}, \qquad \text{with} \quad C_i=[ 1 \ \cdots\ 1] \in \CC^{1\times m_i}.
\end{equation}  
Similarly, let 
$M= [ M_1  \ \cdots \ M_s] \in \ZZ^{n\times m}$
be the block matrix where the columns of $M_i\in \ZZ^{n\times m_i}$ are the elements of $\S_i$ in some fixed order. 
For the vertical system $F=C(a\star x^M)$, it then holds that $\F_\mathrm{free}(\S_1,\dots,\S_s)=\{F_a:a\in \CC^{m}\}$. 
We refer to this $F$ as the \term{vertical system associated with $\S_1,\dots,\S_s$}.

{Let $V_1,\ldots,V_s\subseteq\RR^n$ be sets of vectors. Following the terminology in \cite{Perfect1969,Yu2016prime}, we define an \term{independent transversal} of $(V_1,\ldots,V_s)$ to be a linearly independent tuple $(v_1,\ldots,v_s)\in\prod_{i=1}^s V_i$. It is a well-known linear algebra fact (see, e.g., \cite[Theorem~1]{Perfect1969} and \cite[Theorem~4]{khovanskii2016newton}) that the existence of such an independent transversal is equivalent to 
\begin{equation}\label{eq:essential_family}
\dim\!\Big(\sum_{j\in J} \operatorname{span}_\RR(V_i)\Big)\leq |J|\quad\text{for all $J\subseteq[s]$}.    
\end{equation}
Consider now support sets $\S_1,\ldots,\S_s\subseteq\ZZ^n$, and let $$\Lin(\S_i):=\spn_{\RR}\{u-v:u,v\in\S_i\}\subseteq \RR^n$$ 
denote the 
direction of the affine hull of 
$\S_i$. Then $(\S_1,\ldots,\S_s)$ is called an \term{essential family} in \cite{Sturmfels1994,spaenlehauer-bender} if   $(\Lin(\S_1),\ldots,\Lin(\S_s))$ satisfies \eqref{eq:essential_family}. It is shown in \cite[Lemma~1]{Yu2016prime} and \cite[Theorem~11]{khovanskii2016newton} that this completely characterizes when $\F_\mathrm{free}(\S_1,\ldots,\S_s)$ is generically consistent. We recover this fact as a corollary of \Cref{thm:main_theorem_for_vertical_systems}.
}

\begin{corollary}\label{cor:freely_parametrized}
Let $\S_1,\dots,\S_s\subseteq \ZZ^{n}$, and  let $\mathcal{G}_i\subseteq \CC^n$ be a generating set of $\Lin(\S_i)$ for $i=1,\dots,s$.
The following statements are equivalent:
 \begin{enumerate}[label=(\roman*)]
    \item[(i)] The dimension of $\VV_{\CC^*}(F)$ is $n-s$ for generic $F\in \F_\mathrm{free}(\S_1,\dots,\S_s)$.
    \item[(ii)] {The tuple $(\mathcal{G}_1,\ldots,\mathcal{G}_s)$ admits an independent transversal.}
    \item[(iii)]  {The tuple $(\Lin(\S_1),\ldots,\Lin(\S_s))$ admits an independent transversal.}
\end{enumerate}
Furthermore: 
\begin{itemize}
    \item If the equivalent statements hold, then  $\VV_{\CC^*}(F)$ is  pure-dimensional and all zeros of $F$ are nondegenerate for generic $F\in \F_\mathrm{free}(\S_1,\dots,\S_s)$.
    \item If the equivalent statements do not hold, then $\VV_{\CC^*}(F)$ is empty for generic $F\in \F_\mathrm{free}(\S_1,\dots,\S_s)$, with dimension strictly larger than $n-s$ if nonempty.
    \item The ideal $\langle F\rangle$ is radical for generic $F\in \F_\mathrm{free}(\S_1,\dots,\S_s)$.
\end{itemize}
\end{corollary}

\begin{proof}
Let $F=C(a\star x^M)$ be the vertical system associated with $\S_1,\dots,\S_s$. 
By \Cref{thm:main_theorem_for_vertical_systems}, for $\P=\CC^m$ and $\X=(\CC^*)^n$, statement (i) is equivalent to  
$\rk(C\diag(Gu)M^\top)=s$ for some $u\in \CC^{m-s}$, with $G$ a matrix whose columns form a basis of $\ker(C)$.
By the form of $C$, we can choose $G$ to be
{ \[G = \begin{bmatrix}
    G_1 & 0 & \dots & 0 \\ 0 & G_2 & \dots & 0 \\ \vdots & \vdots & \ddots & 0 \\ 0 & 0 & \dots & G_s
\end{bmatrix},  \qquad \text{with}\quad G_i=\begin{bmatrix}
    1 & 0 & \dots & 0 \\ 0 & 1 & \dots & 0 \\ \vdots & \vdots & \ddots & \vdots \\ 
    0 & 0 & \dots & 1 \\
    -1 & -1 & \dots & -1
\end{bmatrix} \in \CC^{m_i \times (m_i-1)}.\]}
We index $u\in \CC^{m-s}$ in block form as $u=(u_1,\dots,u_s)$ with $u_i=(u_{i1},\dots,u_{i,m_i-1})\in \CC^{m_i-1}$. 
Then, the product
$C\diag(Gu)$ is the block diagonal matrix of the same shape as $C$ but with blocks given by 
\[C_i \diag(G_iu_i)= [u_{i1} \quad \dots \quad u_{i,m_i-1}\quad -(u_{i1}+\dots+u_{i,m_i-1})].\] 
The $i$-th row of $C\diag(Gu)M^\top$ is then the row vector
\[\big( u_{i1}(M_{i1}-M_{i,m_i}) + \cdots +u_{i,m_i-1}(M_{i,m_i-1}-M_{i,m_i}) \big)^\top.   \]
The set of vectors obtained by varying all the entries of $u_i$ is precisely $\Lin(\S_i)$. 
Hence, we have shown that  $\rk(C\diag(Gu)M^\top)=s$ for some $u\in \CC^{m-s}$ (i.e. (i) holds) if and only if (iii) holds. 

Clearly, (ii) implies (iii). For the reverse implication, assume that (ii) does not hold. 
Then, for all choices of $\omega_i\in \mathcal{G}_i$, $i=1,\dots,s$,  the matrix with columns $\omega_1,\dots,\omega_s$ has rank smaller than $s$; that is, all $s$-minors are zero. 
For any $(v_1,\ldots,v_s)\in \prod_{i=1}^s \Lin(\S_i)$, 
we have that $v_i$ is a linear combination of the elements of $\mathcal{G}_i$. Then, the multilinear expansion of the determinant gives 
that each $s$-minor of the matrix with columns $v_1,\dots,v_s$ is a linear combination of $s$-minors of matrices with $i$-th column in  
$\mathcal{G}_i$, $i=1,\dots,s$. Hence it is zero, showing that (iii) does not hold. This concludes the proof of (i) $\Leftrightarrow$ (ii) $\Leftrightarrow$ (iii). 

The bullet points are a direct consequence of \Cref{thm:main_theorem_for_vertical_systems}.
\end{proof}

\begin{example}\label{ex:1}
Consider the freely parametrized system 
\[\widetilde{F}=\left(\begin{array}{l}
 a_{1} x_{1} x_{3} + a_{2} x_{2} x_{3} + a_{3} x_{3} + a_{4}\\
 a_{5} x_{1}^2 + a_{6} x_{2} x_{3} + a_{7} x_{3} + a_{8} \\
 a_{9} x_{1}^2 + a_{10} x_{2} x_{3} + a_{11} x_{3} + a_{12} 
\end{array}\right)\]
with supports
\[\S_1=\{(1,0,1),(0,1,1),(0,0,1),(0,0,0)\},\quad \S_2=\S_3=\{(2,0,0),(0,1,1),(0,0,1),(0,0,0)\}\,.\]
This system is generically consistent since there exists a linearly independent tuple
\[\Big((1,0,1),(2,0,0),(0,1,1)\Big)\in\Lin(\S_1)\times\Lin(\S_2)\times\Lin(\S_3).\]
\end{example}

\begin{remark}
If a vertical system $F$ is generically consistent, so is any other vertical system $\widetilde{F}$ with the same supports and less dependencies between the coefficients (systems \eqref{eq:free_example} and \eqref{eq:vertical_example} in the introduction provide an example of this). However, restricting a vertical system by keeping the supports and adding dependencies among the coefficients does not preserve generic consistency. For example, consider the vertical system
\begin{equation}\label{eq:vertical_inconsistent}
    F = \left(\begin{array}{l} a_{1} x_{1} x_{3} + a_{2} x_{2} x_{3} + a_{3} x_{3} + a_{4}\\
 a_{5} x_{1}^2 + a_{2} x_{2} x_{3} + a_{3} x_{3} + a_{4} \\
 a_{6} x_{1}^2 + a_{2} x_{2} x_{3} + a_{3} x_{3} + a_{4} \end{array}\right),
\end{equation}
which can be seen as a specialization of the generically consistent system from \Cref{ex:1}. 
It is clear from inspection that the zero locus is empty unless $a_5=a_6$, 
and in this case, the dimension of the zero locus is $1$. The system is thus generically inconsistent. 
\end{remark}

The specialization of \Cref{thm:main_theorem_for_vertical_systems} for vertical systems to the case $\A=\RR^m_{>0}$ and $\X=\RR^n_{>0}$ allows us to study freely parametrized systems with real coefficients with prescribed sign. Specifically, 
we consider maximal dimensional   orthants $\mathcal{O}_i$ of $(\RR^*)^{\S_i}$, and consider the  subset $\F_{\mathcal{O}_1,\dots,\mathcal{O}_s}(\S_1,\dots,\S_s)$ of $\F_\mathrm{free}(\S_1,\dots,\S_s)$ consisting of the polynomials in  the image of $\prod_{i=1}^s \mathcal{O}_i$ under the isomorphism $\prod_{i=1}^s\CC^{\S_i}\cong \F_\mathrm{free}(\S_1,\dots,\S_s)$. 

\begin{corollary}\label{cor:fixedsignfree} 
Let $\S_1,\dots,\S_s\subseteq \ZZ^{n}$. For $i=1,\dots,s$, let $\mathcal{G}_i\subseteq \CC^n$ be a generating set of $\Lin(\S_i)$ and  $\mathcal{O}_i$ be a maximal dimensional orthant of $(\RR^*)^{\S_i}$. The following statements are equivalent:
 \begin{enumerate}[label=(\roman*)]
    \item The semialgebraic set $\VV_{\RR^*}(F)\cap \RR^n_{>0}$ has dimension $n-s$ for $F$ in a subset of\\ $\F_{\mathcal{O}_1,\dots,\mathcal{O}_s}(\S_1,\dots,\S_s)$ with nonempty Euclidean interior.
    \item {The tuple $(\mathcal{G}_1,\ldots,\mathcal{G}_s)$ admits an independent transversal} and $\mathcal{O}_i\notin \{ \RR^{\S_i}_{>0},\RR^{\S_i}_{<0}\}$ for all $i=1,\dots,s$.
\end{enumerate}
\end{corollary}
\begin{proof}
Let $F=C(a\star x^M)$ be the vertical system associated with $\S_1,\dots,\S_s$, and let $\widetilde{C}$ be obtained by multiplying by $-1$ each column of $C_i$ in \eqref{eq:Cfree} corresponding to a negative coordinate of $\mathcal{O}_i$. 
Then $\F_{\mathcal{O}_1,\dots,\mathcal{O}_s}(\S_1,\dots,\S_s)$
is identified with the vertical system $\widetilde{C}(a\star x^M)$ with $\A=\RR^m_{>0}$. The result now follows  from \Cref{cor:freely_parametrized}, \Cref{thm:main_theorem_for_vertical_systems} and \Cref{prop:nonempty_intersection_with_incidence_variety}, after realizing that $\ker(\widetilde{C}) \cap \RR^n_{>0}\neq \varnothing$ if and only if $\mathcal{O}_i\neq \RR^{\S_i}_{>0},\RR^{\S_i}_{<0}$ for all $i=1,\dots,s$. 
\end{proof}

In the freely parametrized setting, it is known from \cite[Proposition~1.12]{spaenlehauer-bender} that one can obtain the generic properties predicted by \Cref{cor:freely_parametrized}, by only requiring the coefficients corresponding to the vertices of the Newton polytopes to be generic (while the other coefficients can be arbitrarily fixed). It is an interesting problem for future research to determine which parameters can be fixed in a vertical system while still retaining the generic properties given by \Cref{thm:main_theorem_for_vertical_systems}.

\subsection{Extension to $\CC$}\label{subsec:extension_to_C}
A natural question  is to what extent \Cref{thm:main_theorem_for_vertical_systems} also holds for zeros in the whole affine space $\CC^n$  rather than just in the torus $(\CC^*)^n$. As explained in \Cref{rk:extension_C}, we can extend \Cref{thm:summary_of_implications} to zeros in $\CC^n$, but our proofs of the irreducibility of $\I$ for augmented vertical systems and of \Cref{thm:degenerate_zeros_of_vertical_systems} rely heavily on the fact that we can invert the entries of $x\in (\CC^*)^n$. 
The following example illustrates that the  affine incidence variety 
\[\I_{\CC}=\{(a,b,x)\in\CC^m\times\CC^\ell\times\CC^n:F_{a,b}(x)=0\}\]
might  be reducible for an augmented vertical system. 

\begin{example}\label{ex:non-independent_constant_terms}
For the following variation of Example \eqref{eq:vertical_inconsistent}
\[\left(\begin{array}{l}
 a_{1} x_{1} x_{3} + a_{2} x_{2} x_{3} + a_{3} x_{3} + a_{4}\\[0.1em]
 a_{5} x_{1}^2 + a_{2} x_{2} x_{3} + a_{3} x_{3} + a_{4} \\[0.1em]
 a_{6} x_{1}^3 + a_{2} x_{2} x_{3} + a_{3} x_{3} + a_{4} 
\end{array}\right),\]
the affine incidence variety has two irreducible components:
\[
\I_{\CC}=\VV_{\CC}(x_{1},
a_{2} x_{2} x_{3} + a_{3} x_{3} + a_{4})\cup 
\VV_{\CC}( a_{1} x_{3} - a_{6} x_{1}^2, a_1 x_3 - a_5 x_1, 
  a_{1} x_{1} x_{3} + a_{2} x_{2} x_{3} + a_{3} x_{3} + a_{4}).\] 
The system is generically consistent over $\CC^*$ (and thus also over $\CC$), and hence the zero locus over $\CC^*$ is generically $0$-dimensional. However, the zero locus over $\CC$ is  generically 1-dimensional, due to the curve  in the $\{x_1=0\}$ coordinate hyperplane, which contradicts the expected dimension of $n-s=0$. 
\end{example}

A simple criterion that guarantees that the affine incidence variety  is irreducible   is that each polynomial appearing in $F$ has a constant term involving a unique parameter. More precisely, we have the following result, where (i) is an affine analog of \Cref{thm:irreducible_smooth_incidence_variety}, 
(ii) is an affine analog of \Cref{thm:main_theorem_for_vertical_systems}, and (iii) can be seen as an extension of \cite[Lemma~2.1]{LW96} to augmented vertical systems. 

\begin{theorem}\label{thm:extension_to_C}
Let $F=(C(a\star x^M),Lx-b)\in \CC[a,b,x]^{s+\ell}$ be an augmented vertical system  with $C\in \CC^{s\times m}$ of full rank $s$, $L\in \CC^{\ell\times n}$ with $s+\ell\leq n$ and $M\in \ZZ_{\geq 0}^{n\times m}$. Suppose 
that for some indices $i_1<\dots<i_s$, the submatrix of $C$ given by the columns with these indices is diagonal of rank $s$, and that the corresponding submatrix of $M$ is the zero matrix. Then the following holds:
 \begin{enumerate}[label=(\roman*)]
\item The affine incidence variety $\I_{\CC}$ is nonsingular and irreducible of dimension $m+n-s$. 
\item Suppose that $(\A\times\B,\X)\subseteq \CC^m\times \CC^\ell\times \CC^n$ is an algebraically defined pair, and assume that 
$\I_{\CC}\cap (\A\times\B\times \X)$ is Zariski dense in $\I_{\CC}$. Then the affine analogs of {\rm (deg1)}, {\rm (degX1)}, {\rm (degXG)}, {\rm (degAllG)}, {\rm (setE)}, {\rm (setZ)}, {{\rm (flatG)}}, {\rm (dim1)}, {\rm (dimG)} are equivalent, any of these statements imply {\rm (dimX)}, {\rm (real)} and {\rm (reg)}, and {\rm (rad)} holds independently of the other statements. 
\item For generic $(a,b)\in\CC^m\times\CC^\ell$, the variety $\VV_{\CC}(F_{a,b})$ has no irreducible component contained in a coordinate hyperplane of $\CC^n$.
\end{enumerate}
\end{theorem}

\begin{proof}
To prove (i), we parametrize the affine incidence variety as follows.
Without loss of generality, we  assume that the first $s$ columns of $C$ form the diagonal matrix. Then $a_1,\dots,a_s$ are terms of the entries $1,\dots,s$ of $C(a \star x^M)$ respectively. Hence  $C(a \star x^M)=0$ is equivalent to
$(a_1,\dots,a_s)=\psi(a',x)$ where $a'=(a_{s+1},\dots,a_m)$ and $\psi$ is a polynomial function. This gives 
the parametrization $(\psi(a',x),a',Lx, x)$
of $\I_{\CC}$  in terms of $m+n-s$ free parameters.  
To show nonsingularity, we proceed analogously to the proof of \Cref{thm:irreducible_smooth_incidence_variety}(ii), after noting that the first $s$ columns of $C \diag(x^M)$ are independent of $x\in \CC^n$ and form a diagonal matrix of full rank $s$.

For statement (ii), the equivalences in \Cref{thm:summary_of_implications} and the equivalence between (setZ) and (setE) hold by part (i),  \Cref{rk:extension_C} and the extension of \Cref{prop:properties_of_alg_defined_pairs}(ii) to $\CC$.  
The only obstruction lies in proving that  (setZ) implies (degAllG); with that in place, the rest of the claims in (ii)  follow analogously to the proof of \Cref{thm:main_theorem_for_vertical_systems}. 

The proof of (setZ) $\Rightarrow$(degAllG)  is based on the following construction.
For every  subset $I\subseteq \{1,\dots,n\}$, let $F_I\in \CC[a,b,x]^{s+\ell}$ be the system obtained by letting $x_i=0$ for $i\notin I$ in $F$. As each entry  of $F_I$ has a constant term, the coefficient matrix of the vertically parametrized part of $F_I$ has maximal rank. Hence $F_I$  is again an augmented vertical system. 

Let $\X^*:=\X\cap (\CC^*)^n$ and let $\Theta_I\subseteq \CC^n$ be the set defined by $x_i=0$ if $i\notin I$. 
For $(a,b)\in \CC^m\times\CC^\ell$, there are surjective maps $\psi_I\from \VV_{\CC^*}(F_{I,a,b}) \to \VV_{\CC}(F_{a,b}) \cap \Theta_I$
obtained by sending $x^*\in \VV_{\CC^*}(F_{I,a,b})$ to $x'\in \Theta_I$ defined as $x'_i=x^*_i$ if $i\in I$ and zero otherwise. 
We obtain in particular that
\begin{equation}\label{eq:decomp}
\VV_{\CC}(F_{a,b}) =  \bigcup_{ I\subseteq \{1,\dots,n\}} \psi_I(\VV_{\CC^*}(F_{I,a,b})) \quad\text{and}\quad 
\Z_{F,\P}(\X) =\bigcup_{I\subseteq \{1,\dots,n\}} \Z_{F_I,\P}(\X^*).
\end{equation}
(Here and below, subindices $F,F_I$ are added to differentiate among the augmented vertical systems.)
Additionally, 
the columns of $J_{F_{a,b}}(\psi_I(x^*))$ and of $J_{F_{I,a,b}}(x^*)$ indexed by $I$ agree, and $J_{F_{I,a,b}}(x^*)$  is zero outside these columns. Hence, 
if $x^*$ is a nondegenerate  zero of $F_{I,a,b}$, then so is $\psi_I(x^*)$  as a zero of $F_{a,b}$. 

Set $\P=\A\times \B$. If (setZ) holds for $F$,  then $\Z_{F_I,\P}(\X^*)$ is Zariski dense for at least one $F_I$ by the second equality in \eqref{eq:decomp}. In this case, so is $\Z_{F_I}=\Z_{F_I,\CC^m\times \CC^\ell}((\CC^*)^n)$ and hence, using the implication 
(setZ) $\Rightarrow$(degAllG) from \Cref{thm:main_theorem_for_vertical_systems} for $F_I$, there exists a nonempty Zariski open subset $U_I\subseteq \CC^m\times\CC^\ell$
such that all zeros of $F_{I,a,b}$ in $(\CC^*)^n$ are nondegenerate for $(a,b)\in U_I$. 
 Consider the nonempty Zariski open set $U$   obtained by intersecting $U_I$ over all subsets $I$ for which $\Z_{F_I,\P}(\X^*)$ is Zariski dense, and removing  the Zariski closed sets 
$\overline{\Z_{F_I,\P}(\X^*)}$ for all other subsets $I$. 
Then, for all $(a,b)\in U$ and for all $I$, all zeros of 
$F_{I,a,b}$ in $(\CC^*)^n$ are nondegenerate, and hence so are all zeros of $F_{a,b}$ in $\CC^n$, which gives (degAllG). 

Finally, to prove (iii), we begin by noting that the statement is trivial if $F$ is generically inconsistent. Hence, we assume $F$ is generically consistent. If we set $x_i=0$, we obtain an augmented vertically parametrized system $F_{\vert x_i=0}$ with independent constant terms and $n-1$ variables. 
Suppose this system is generically consistent. Then, for generic $(a,b)$, all irreducible components of $\VV_{\CC}(F_{a,b})\cap\{x_i=0\}\cong \VV_{\CC}((F_{\vert x_i=0})_{a,b})$ have dimension $n-1-(s+\ell)$ by part (ii). But by (ii) again, all components of $\VV_{\CC}(F_{a,b})$ have dimension $n-(s+\ell)$. Hence, none of them is fully contained in $\{x_i=0\}$ for generic  $(a,b)$.
\end{proof}

In the freely parametrized setting, we recover  \cite[Proposition~2.1]{Labahn2021homotopy}.

\begin{corollary}
For $\S_1,\ldots,\S_s\subseteq\ZZ_{\geq0}^n$ with $0\in \S_i$ for $i=1,\dots,s$, the following are equivalent:
\begin{enumerate}[label=(\roman*)]
\item $\VV_{\CC}(F)$ has pure dimension $n-s$ for generic $F\in\F_\mathrm{free}(\S_1,\ldots,\S_s)$.
\item $\VV_{\CC^*}(F)$ has pure dimension $n-s$ for generic $F\in\F_\mathrm{free}(\S_1,\ldots,\S_s)$. 
\item {The tuple $(\S_1,\ldots,\S_s)$ admits an independent transversal.}
\end{enumerate}
\end{corollary}
\begin{proof}
Note that $\S_i$ is a generating set of $\Lin(\S_i)$ as $0\in \S_i$. Hence (ii) and (iii) are equivalent by \Cref{cor:freely_parametrized}. 

Let $H$ be the vertical system associated with $\S_1,\ldots,\S_s$. 
To show that (i) implies (ii), assume (dimG) holds for $H$ over $\CC$.
If $m$ is the sum of the cardinalities of $\S_1,\ldots,\S_s$, taking Zariski closures in affine space, we have 
\[ \overline{\Z((\CC^*)^n)}=\overline{\pi(\I)} =  \overline{\pi(\overline{\I})}=
\overline{\pi(\I_\CC)}=\overline{\Z(\CC^n)} = \CC^m,
\]
where in the third equality we use that $\I_{\CC}$ is irreducible by \Cref{thm:extension_to_C}(i) and that $\I\neq \varnothing$, and in the last equality, we 
use  the implication (dimG)$\Rightarrow$(setZ) over $\CC$ of \Cref{thm:extension_to_C}(ii). So, by \Cref{thm:main_theorem_for_vertical_systems}, $\VV_{\CC^*}(F)$ has pure dimension $n-s$ for generic $F\in\F_\mathrm{free}(\S_1,\ldots,\S_s)$ and hence (ii) holds.

Conversely, if (ii) holds, then (deg1) holds for $H$ over $\CC^*$ by \Cref{cor:freely_parametrized}. Hence, (deg1) holds   for $H$ over $\CC$ as well, and the implication (deg1)$\Rightarrow$(dimG)  over $\CC$ of \Cref{thm:extension_to_C}(ii) gives that $\VV_{\CC}(F)$  has pure dimension $n-s$ for generic $F\in\F_\mathrm{free}(\S_1,\ldots,\S_s)$, that is, (i) holds.  
\end{proof}

\begin{example}\label{ex:independet_constant_terms}  
We modify the generically consistent system from \Cref{ex:non-independent_constant_terms} to
\[\left(\begin{array}{l}
 a_{1} x_{1} x_{3} + a_{2} x_{2} x_{3} + a_{3} x_{3} + a_{4}\\[0.1em]
 a_{5} x_{1}^2 + a_{2} x_{2} x_{3} + a_{3} x_{3} + a_{6} \\[0.1em]
 a_{7} x_{1}^3 + a_{2} x_{2} x_{3} + a_{3} x_{3} + a_{8} 
\end{array}\right),\]
where each polynomial now contains an independent constant term ($a_4$, $a_6$ and $a_8$, respectively). By  \Cref{thm:extension_to_C}(i), the affine incidence variety in $\CC^{8}\times\CC^3$ is irreducible of dimension $8$. 
The system still falls in the generically consistent scenario, and it therefore follows from \Cref{thm:extension_to_C}(ii) that the zero locus over $\CC$ generically has the expected dimension $0$.
\end{example}

An interesting problem for future research is to explore conditions other than those in \Cref{thm:extension_to_C} that ensure that the zero locus of an augmented vertical system generically has no component contained in the coordinate hyperplanes, similar to what has been done in the square freely parametrized setting (see, e.g., \cite{LW96,RW96}).

\appendix
\section{Nondegeneracy and radicality}\label{app:radical}
In this appendix we prove \Cref{prop:generic_radicality_in_Laurent_polynomial_ring}. The proof is based on a series of lemmas, the first of which  is a well-known commutative algebra fact adapted to the Laurent polynomial setting.

\begin{lemma}
\label{prop:nondegeneracy_implies_radicality}
\leavevmode
\begin{enumerate}[label=(\roman*)]
\item  If all zeros of $F\in \CC[x]^r$ in $\CC^n$ are nondegenerate, then $\langle F\rangle\subseteq\CC[x]$ is radical.
\item  If all zeros of $F\in \CC[x^\pm]^r$ in $(\CC^*)^n$ are nondegenerate, then $\langle F\rangle\subseteq\CC[x^\pm]$ is radical.
\end{enumerate}
\end{lemma}
\begin{proof}
Part (i) follows by \cite[Corollary~16.20]{eisenbud1995commutative} (which, in turn, is a consequence of the principal ideal theorem \cite[Theorem~10.2]{eisenbud1995commutative} and the Jacobian criterion \cite[Theorem~16.19]{eisenbud1995commutative}), combined with the fact that being reduced is a local property of rings.

For (ii), first of all note that multiplying the polynomials $F=(f_1,\ldots,f_r)$ by monomials neither changes the ideal $\langle F\rangle\subseteq\CC[x^{\pm}]$ nor the nondegeneracy of their zero locus over $\CC^*$. Hence,  without loss of generality, we can assume $F \in \CC[x]^r$.
In order to show that $\langle F\rangle\subseteq\CC[x^\pm]$ is radical, it is now enough to show that 
\[I:=\langle f_1,\ldots,f_r,x_1y_1-1,\ldots,x_ny_n-1 \rangle\subseteq\CC[x_1,\ldots,x_n,y_1,\ldots,y_n]=:\CC[x,y]\]
is radical.
The Jacobian matrix of the generators defining $I$ has the block form
\[J:= \begin{bmatrix}
    J_F(x)&0\\
    \diag(y)&\diag(x)
\end{bmatrix} \in \CC^{(r+n) \times 2n},\]
and has rank $r+n$ for all $(x,y)\in\VV_{\CC}(I)$, as for any such point, $x\in (\CC^*)^n$  is a zero of $F$, and hence nondegenerate.
Hence, it follows by (i) that $I$ is radical.
\end{proof}

The following lemma allows us to switch back and forth between radicality in the Laurent polynomial ring and the usual polynomial ring.

\begin{lemma}
\label{lem:radicality_for_Laurant_polynomials_vs_radicality_of_saturation} Let $K$ be a field and let 
$F\in K[x^\pm]^r$.  
Then $\langle F \rangle\subseteq K[x^\pm]$ is radical if and only if $\langle F \rangle\cap K[x]\subseteq K[x]$ is radical.
\end{lemma}

\begin{proof}
The ``only if'' direction holds as contractions of radical ideals are radical. For the ``if'' direction, assume $\langle F\rangle \cap K[x]$ is radical and that $f^u\in\langle F\rangle$ for some $f\in K[x^\pm]$ and some integer $u>0$. Let $N>0$ be an integer such that $(x_1\cdots x_n)^{Nu}f^u\in \langle F\rangle\cap K[x]$. Then, by assumption, $(x_1\cdots x_n)^Nf\in\langle F\rangle \cap K[x]$, from which it follows that $f\in \langle F\rangle$.
\end{proof}

For a parametric system $F\in \CC[p,x^\pm]^r$ with $p=(p_1,\dots,p_k)$, we next relate  radicality for generic parameter values $p$ to radicality over the field $\CC(p)$ of rational functions in the  parameters.

\begin{lemma}
\label{lem:generic_radicality_vs_radicality_over_function_field_saturation}
Let $F=(f_1,\dots,f_r)\in \CC[p,x^\pm]^r$. If $\langle F_p\rangle\subseteq\CC[x^\pm]$ is radical for generic $p\in\CC^k$, then $\langle F \rangle\subseteq\CC(p)[x^\pm]$ is a radical ideal.
\end{lemma}

\begin{proof}
The proof will rely on Gröbner bases, so we will reformulate the lemma from a statement about Laurent polynomial rings to a statement about usual polynomial rings. 
We will use a subscript $R$ in the ideal notation $\langle\cdot\rangle_R$ to indicate that the ideal is generated in the ring $R$ whenever this is not clear from the context.
Since monomials in $x$ are units in $\CC[p,x^\pm]$, we can without loss of generality assume that  $F\in \CC[p,x]^r$. Then contraction corresponds to saturation: 
\[\begin{array}{l}\langle F\rangle_{\CC(p)[x^\pm]}\cap\CC(p)[x]=\langle F\rangle_{\CC(p)[x]}{:}(x_1\cdots x_n)^\infty,\\[0.5em]
\langle F_p\rangle_{\CC[x^\pm]}\cap\CC[x]=\langle F\rangle_{\CC[x]}{:}(x_1\cdots x_n)^\infty\:\:\text{for $p\in\CC^k$}.
\end{array}\]
By \Cref{lem:radicality_for_Laurant_polynomials_vs_radicality_of_saturation}, the lemma at hand says that if $I_p :=\langle F_p\rangle_{\CC[x]}{:}(x_1\cdots x_n)^\infty$ is radical for generic $p\in\CC^k$, then  $I:=\langle F \rangle_{\CC(p)[x]}{:}(x_1\cdots x_n)^\infty$ is radical.

We can  construct a generating set of $I$ that generically specializes to 
a generating set of $I_p$ in the following way.
Let $\widetilde{G}$ be  a Gröbner basis of  $\langle f_1,\dots,f_r, 1-x_1\cdots x_n\, y\rangle$ in the ring $\CC(p)[x_1,\ldots,x_n,y]$ with respect to the lexicographic ordering $y>x_1 >\cdots > x_n$. 
Then $G:=\widetilde{G}\cap \CC(p)[x]$ is a Gröbner basis for $I$ by \cite[Theorem~4.4.14]{cox2015ideals}. Now $\widetilde{G}$ specializes to a Gröbner basis $\widetilde{G}_{p}$ of $\langle f_{1,p},\dots,f_{r,p}, 1-x_1\cdots x_n\, y \rangle $ for generic $p$ by \cite[Theorem~6.3.1]{cox2015ideals}, so that $\widetilde{G}_{p}\cap \CC[x]$ is a Gröbner basis for $I_p$. Also, $\widetilde{G}_{p}\cap \CC[x]=G_{p}$ for generic $p$. Hence, there exists a nonempty Zariski open set $U\subseteq \CC^{k}$ such that 
$G$ specializes to a Gröbner basis $G_{p}$ of $I_p$ for all $p\in U$. 

We now prove the contraposition of the desired result. If $I$ is not radical, then there exists   $h\in \CC(p)[x]$ 
and  $N \geq 0$ such that 
the normal form of $h^N$ with respect to $G$ is zero, while the normal form of $h$ is nonzero. 
Let $Z$ be the proper Zariski closed subset of $\CC^k$ where the denominators of the quotients and remainders of the division of $h^N$ and $h$  by $G$ vanish, and  where the normal form of $h$ vanishes. Then, for all $p$ in the nonempty  Zariski open subset $U\setminus Z$, the normal form of $h_p$ by $G_p$ is nonzero and that of $h_p^N$ is zero. Hence, $h_p\notin I_p$ but $h_p\in \rad(I_p)$. This gives a contradiction, showing the statement.
\end{proof}

\begin{proof}[{Proof of \Cref{prop:generic_radicality_in_Laurent_polynomial_ring}}]
If all zeros of $F_p$ are nondegenerate for generic $p\in \CC^{k}$, then, by \Cref{prop:nondegeneracy_implies_radicality}, we have that $\langle F_p\rangle\subseteq\CC[x^{\pm}]$ is  radical for generic $p\in \CC^{k}$. 
From this, \Cref{lem:radicality_for_Laurant_polynomials_vs_radicality_of_saturation}, gives that $\langle F_p\rangle\cap \CC[x]\subseteq\CC[x]$ also is generically radical in $\CC[x]$, and by \Cref{lem:generic_radicality_vs_radicality_over_function_field_saturation}, the ideal  $\langle F\rangle\subseteq\CC(p)[x^\pm]$ is radical. Using again \Cref{lem:radicality_for_Laurant_polynomials_vs_radicality_of_saturation} with the field $\CC(p)$, we obtain that $\langle F\rangle\cap\CC(p)[x]\subseteq\CC(p)[x]$ is radical.
\end{proof}

\newlength{\bibitemsep}\setlength{\bibitemsep}{.2\baselineskip plus .05\baselineskip minus .05\baselineskip}
\newlength{\bibparskip}\setlength{\bibparskip}{0pt}
\let\oldthebibliography\thebibliography
\renewcommand\thebibliography[1]{
  \oldthebibliography{#1}
  \setlength{\parskip}{\bibitemsep}
  \setlength{\itemsep}{\bibparskip}
}

\bigskip

\bibliographystyle{alpha} 

\newcommand{\etalchar}[1]{$^{#1}$}

\bigskip

\small
\noindent {\bf Authors' addresses:}

\noindent 
Elisenda Feliu, University of Copenhagen \hfill{\tt efeliu@math.ku.dk}\\
Oskar Henriksson, University of Copenhagen \hfill{\tt oskar.henriksson@math.ku.dk}\\
Beatriz Pascual-Escudero, Universidad Politécnica de Madrid \hfill {\tt beatriz.pascual@upm.es}

\end{document}